\input amssym.def
\input amssym
\documentclass{tac}

\usepackage{xcolor}

\title{On Deformations of Pasting Diagrams}

\copyrightyear{2009}

\thanks{The author wishes to thank Kansas State University and the 
University of Pennsylvania for support of the 
sabbatical leave when this work was begun.  The author also wishes to 
than the referee assigned by {\em Theory and Applications of 
Categories} for numerous corrections and helpful suggestions for 
improvements to the original manuscript.}

\eaddress{dyetter@math.ksu.edu}

\author{D.\ N.\ Yetter} 

\address{Department of Mathematics \\
Kansas State University \\ Manhattan, KS 66506}

\begin{document}
\bibliographystyle{plain}

\maketitle

\begin{abstract}  We adapt the work of Power \cite{Pow} to
describe general, not-necessarily composable, not-necessarily
commutative 2-categorical pasting diagrams and their composable and
commutative parts.  We provide a deformation theory for 
pasting diagrams valued in the 2-category of $k$-linear categories, 
paralleling that
provided for diagrams of algebras by Gerstenhaber and Schack
\cite{GSdiag}, proving the standard results.  Along the
way, the construction gives rise to a bicategorical analog of the
homotopy G-algebras of Gerstenhaber and Voronov \cite{GV}.
\end{abstract}

\tableofcontents

\section{Introduction}

It is the purpose of this work to describe the deformation theory
of pasting diagrams of $k$-linear categories, functors and natural
transformations.  As such, it generalizes work of Gerstenhaber and
Schack \cite{GSdiag}, both by giving an exposition of the well-known
extension
of their work on diagrams of algebras to $k$-linear categories and
functors, and by the inclusion of natural transformations.

The present work has a number of motivations.  It 
initially grew out of a program to extend the author's 
deformation theory for monoidal categories, functors and
natural transformations \cite{CY.def, Y.def, Y.book, Y.abkan},
which deforms only the structure maps, to a theory in which
the composition of the category, the arrow part of the
monoidal product, and the structure maps are all deformed
simultaneously.  That program is still in progress, and
this paper is a first step in it.

It is also the first step in a program to provide a 
Gerstenhaber-style deformation theory for linear stacks, as
pre-stacks may be considered as special instances of pasting
diagrams.  Consideration of not-necessarily abelian linear stacks is
motivated by physical considerations in a prospective
deformation quantization approach to quantum gravity (cf.
\cite{BB, Cr1, Cr2}).

It is also hoped that the present work may shed light, if only
by analogy, on the difficulties arising in Elgueta's deformation
theory for monoidal bicategories \cite{Elg}.

Throughout we will consider all categories to be small, if necessary
by invoking the axiom of universe.  Composition will be written in
diagrammatic order unless parentheses indicate
functional application.  Thus $fg$ and $g(f)$ both denote the arrow
obtained by following $f$ with $g$, in the second case thought of as
applying $g$ to $f$.  Throughout by abuse of notation, whenever
operations defined as applying to 
sets, are applied to elements, the singleton set will be understood.
$k$ will denote a fixed field, and all categories and functors will 
be $k$-linear.

\section{Pasting schemes and pasting diagrams: definitions}

A pasting diagram is to $n$-categories what an ordinary diagram is
to categories.  A number of ways to formalize them have been
developed.
We will for the most part follow Power \cite{Pow}, whose approach
mixing Street's notion of
computads \cite{street87} with a geometric adaptation of Johnson's 
pasting schemes \cite{johnson88pasting} avoids much of the
combinatorial complexity
of Johnson's approach. Power's description seems to be of the right
generality for the 
present work, and we will deviate from it only to allow the
description of not-necessarily composable, not-necessarily
commutative diagrams. We follow Power's gentle method of exposition
by initially restating the familiar
in less-familiar but readily generalizable terms:

\begin{definition}
A {\em $1$-computad} is a (finite) directed graph .  A {\em
$1$-computad
morphism} is a map of directed graphs.
\end{definition}

Observe, in particular, that there is a forgetful functor from the
category of small (or finite) categories and functors to
$1$-computads. 

\begin{definition}
A {\em $1$-pasting scheme} is a finite non-empty set $G$ equipped
with an
embedding to the oriented line $\Bbb R$.  The elements (identified
with their
images in the line) are called {\em $0$-cells} of the pasting scheme,
and the open bounded 
intervals in ${\Bbb R}\setminus G$ are called the {\em $1$-cells} of
the
pasting scheme.

Denoting the sets of $0$- and $1$-cells by $G_0$ and $G_1$
respectively, there is a function $dom:G_1\rightarrow G_0$, (resp.
$cod:G_1\rightarrow G_0$) the {\em domain} (resp.
{\em codomain} function, which assign the
lesser (resp. greater) endpoint to each $1-cell$.  

The quadruple $(G_0,G_1,dom,cod)$ defines a finite directed graph,
{\em the underlying 1-computad} of the 1-pasting scheme, which we
denoted $C(G)$.

The {\em domain of G}, $dom G$ (resp. {\em codomain of G} $cod G$) is
the
smallest (resp. largest) element of $G_0$.  (These should not be
confused with
the domain and codomain functions.)
\end{definition}

\begin{definition}
A {\em composable $1$-pasting diagram} in a $1$-computad $H$ consists
of a $1$-pasting scheme and a $1$-computad morphism
$h:C(G)\rightarrow H$.  The {\em domain}
(resp. {\em codomain}) of the $1$-pasting diagram is $h(dom G)$ (resp.
$h(cod G)$).  In the case where $H$ is the underlying $1$-computad of
a category $\cal A$, we call the $1$-pasting diagram a {\em labelling
of
$G$ in $\cal A$}.  We denoted the set of composable $1$-pasting
diagrams in a $1$-computad $H$ by $diags(H)$.
\end{definition}

\begin{definition}
The {\em $1$-pasting composition} of a labelling of a pasting scheme
$G$ in $\cal A$ is the arrow of $\cal A$ obtained by composing the
sequence of arrows in $A$ given by the $h(G_1)$.
\end{definition}

\begin{definition}
A {\em general $1$-pasting diagram} in a $1$-computad $H$ consists of
a $1$-computad
$G$ and a $1$-computad morphism $h:G\rightarrow H$.  The {\em
composable
parts} of a general $1$-pasting diagram $(G,h)$ are all composable
$1$-pasting
diagrams $p:C(S)\rightarrow H$ in $H$ that factor as
$C(S)\stackrel{\tilde{p}}{\rightarrow} G \stackrel{h}{\rightarrow} H$.
\end{definition}

Thus far we have not really added anything to the notions of iterated
composition and diagrams in a category, save to emphasize a sense
that it is somehow
`more geometric' than multiplication.  This heretofore needless
abstraction now becomes necessary:

\begin{definition} A {\em $2$-computad} consists of a $1$-computad $G$ and
a set $G_2$ together with functions $dom:G_2\rightarrow diags(G)$
(resp. $cod:G_2\rightarrow diags(G)$) such that $dom(dom) = dom(cod)$
and $cod(dom) = cod(cod)$ as functions from $G_2$ to $G_0$.

A {\em $2$-computad morphism} from $(G,G_2,dom,cod)$ to
$(H,H_2,dom,cod)$
consists of a $1$-computad morphism $f:G\rightarrow H$ and a function
$f_2:G_2\rightarrow H_2$ such that $dom(f_2) = f(dom)$ and $cod(f_2)
= f(cod)$.
\end{definition}

Plainly $2$-computads and $2$-computad morphisms form a category, and
every
$2$-category has an underlying $2$-computad obtained by forgetting the
compositions.

\begin{definition} A {\em $2$-pasting scheme} consists of a finite
$1$-computad
$G$ together with an embedding of the underlying geometric directed
graph (also denoted $G$ by abuse of notation)
into the oriented plane ${\Bbb R}^2$ satisfying the following 
conditions:

\begin{enumerate}
\item The complement of the image of $G$ consists of an unbounded
region and 
finitely many open cells, which are called {\em faces};
\item The boundary of each bounded face $F$ is of the form $\sigma(F)
\cup
-\tau(F)$, where the negation indicates orientation reversal, each of
$\sigma(F)$ and $\tau(F)$ are images of composable $1$-pasting
diagrams 
in $G$, and $dom \sigma(F) = dom \tau(F)$ (resp. $cod \sigma(F) = cod
\tau(F)$); and
\item There exist vertices $s(G)$ and $t(G)$ in the boundary of the
unbounded face such that for every vertex $v$, there is a composable
$1$-pasting
diagram $h:C(H)\rightarrow G$ such that $v$ is in the image of $h$,
$dom h = s(G)$ and $cod h = t(G)$.
\end{enumerate}

It follows from these conditions that the boundary of the unbounded
face $E$ is also a union of images of composable $1$-pasting diagrams
$-\sigma(E) \cup \tau(E)$ (note the orientation reversal) with
$dom \sigma(E) = dom \tau(E) = s(G)$ and $cod \sigma(E) = cod \tau(E)
= t(G)$.  We define the {\em domain} (resp. {\em codomain}) of 
$G$ by $dom G = \sigma(E)$ (resp. $cod G = \tau(E)$). 
\end{definition}

Every $2$-pasting scheme $G$ admits an underlying $2$-computad $C(G)$
in which
$G$ is the underlying $1$-computad, and $G_2$ is the set of bounded
faces, which are called {\em $2$-cells}.  The domain (resp. codomain)
of a $2$-cell is given by $dom F = \sigma(F)$ (resp. $cod F =
\tau(F)$),
thereby defining the maps $dom, cod:G_2\rightarrow diags(G)$.  Again
these
maps should not be confused with the domain and codomain of the entire
pasting scheme.

\begin{definition}
A {\em composable 2-pasting diagram} in a 2-computad $H$ consists of a
$2$-pasting scheme and a 2-computad morphism $h:C(G)\rightarrow H$.
In the case where $H$ is the underlying $2$-computad of a 2-category
$\cal A$, we call the $2$-pasting diagram a {\em labelling of
$G$ in $\cal A$}.  We denoted the set of composable $2$-pasting
diagrams in a $2$-computad $H$ by $diags(H)$.
\end{definition}

In \cite{Pow} Power proved

\begin{theorem}
Every labelling of a composable 2-pasting scheme has a unique
composite.
\end{theorem}

Meaning that every iterated application of the 1- and 2-dimensional
compositions to the natural transformations in a composable 2-pasting
scheme which results in a single natural transformation gives the
same result.

We will also need

\begin{definition} A {\em general 2-pasting diagram} in a 2-computad $H$
consists of a 2-computad $G$ and a 2-computad morphism
$h:G\rightarrow H$.
The {\em composable
parts} of a general $2$-pasting diagram $(G,h)$ are all composable
$2$-pasting
diagrams in $H$ $p:C(S)\rightarrow H$ that factor as
$C(S)\stackrel{\tilde{p}}{\rightarrow} G \stackrel{h}{\rightarrow} H$.
\end{definition} 

We will not need Power's further explication of corresponding
structures to describe composable pasting diagrams in higher
dimensions.  For our purposes,
it suffices to make

\begin{definition}
A {\em $3$-computad} consists of a $2$-computad $G$ and
as set $G_3$ together with functions $dom:G_3\rightarrow diags(G)$
(resp. $cod:G_3\rightarrow diags(G)$) such that $dom(dom) = dom(cod)$
and $cod(dom) = cod(cod)$ as functions from $G_3$ to
$diags(G_1,G_0,cod,dom)$.

A {\em $3$-computad morphism} from $(G,G_3,dom,cod)$ to
$(H,H_3,dom,cod)$
consists of a $2$-computad morphism $f:G\rightarrow H$ and a function
$f_3:G_2\rightarrow H_2$ such that $dom(f_3) = f(dom)$ and $cod(f_3)
= f(cod)$.
\end{definition}

\noindent and to observe that any $3$-category admits an underlying
$3$-computad.

Of course, a 2-category can be regarded as a 3-category in which all
3-arrows are identities.  By adopting this view, we can use general
3-pasting diagrams, targetted in the underlying 3-computad of a
2-category
to specify commutativity conditions in general 2-pasting diagrams,
since
the presence of a 3-cell asserts the equality of its source and
target.

Thus we make

\begin{definition}
A {\em $k$-linear pasting diagram} is a $3$-computad $G$ together with
a $3$-computad morphism to the underlying $3$-computad of $k-Cat$, the
2-category of all small $k$-linear categories, $k$-linear functors,
and 
natural transformations.
\end{definition}

These are our primary objects of study.
Note that we do not specify the dimension here:  as part of an
abstract hierarchy, these are 3-dimensional objects, but since we are
working in 2-categories, they are degenerate, and in some sense still
2-dimensional. 

\section{Deformations of categories, functors\\ and natural
transformations: definitions and elementary results}

The generalization of Gerstenhaber's deformation theory from
associative algebras \cite{G1, G2} to linear categories, or
'algebroids' in the sense of Mitchell \cite{Mit}, is quite
straight-forward, and both the one readily
available source in which the construction has appeared 
\cite{Anel} and unpublished lectures of Tsygan \cite{Tsy} treat it as
a folk-theorem.

The deformation theory for linear functors, or for that matter
commutative diagrams of linear functors, is similarly a
straight-forward generalization of work of Gerstenhaber and Schack
\cite{GSdiag}.

Finally, the deformation of natural transformations between
undeformed functors is completely trivial, as will be seen.  It is
only when all three elements are combined that there is really
anything new.

To fix notation, we review the basic elements of the theory:

\begin{definition} A {\em deformation $\hat{\cal C}$ of a $k$-linear
category $\cal C$} is an $R$-linear category (i.e. a category 
enriched in $R$-modules), for $R$ a unital 
commutative
local
$k$-algebra, with maximal ideal $\frak m$, whose objects are those of
$\cal C$, with $\hat{\cal C}(X,Y) = {\cal C}(X,Y)\otimes_k R$,
whose composition and identity arrows reduce
modulo $\frak m$ to those of $\cal C$.

For $R$ as above, an {\em $\frak m$-adic deformation $\hat{\cal C}$ of a $k$-linear
category $\cal C$} is a category enriched in the category of
$\frak m$-adically 
complete $R$-modules with obvious the monoidal structure given 
$\hat{\otimes}_{R}$,
the $\frak m$-adic completion of $\otimes_{R}$, whose objects
are those of $\cal C$, with 
$\hat{\cal C}(X,Y) = {\cal C}(X,Y)\hat{\otimes}_k R$, whose 
composition and identity arrows reduce
modulo $\frak m$ to those of $\cal C$.  

Two deformations in either sense are {\em
equivalent} if there is an isomorphism of categories between
them that reduces to the identity functor modulo $\frak m$.  We refer
to such 
an isomorphism as {\em an equivalence of deformations}.

{\em The trivial deformation} is the deformation in which the
composition on the original category is simply extended by
bilinearity, or by bilinearity and continuity in the $\frak m$-adic 
case, while {\em a trivial deformation} is one equivalent to
the trivial deformation.
\end{definition}

Throughout we will be concerned only with {\em $n^{th}$-order
deformations} in which $R = k[\epsilon]/\langle \epsilon^{n+1}
\rangle$ 
and {\em formal deformations}, $\langle \epsilon \rangle$-adically
complete deformations with respect to $R = k[[\epsilon]]$.
Collectively, we
will refer to these as {\em 1-parameter deformations}.

In these cases the composition in $\hat{\cal C}$ has the form

\[ f\star g = \sum \mu^{(j)}(f,g)\epsilon^j \]

\noindent while identity maps in $\hat{\cal C}$ are of the form

\[ \hat{1}_{A} = \sum \iota^{(j)}(A)\epsilon^{j} \]

\noindent where $\mu^{(0)}(f,g) = fg$,  and $\iota^{(0)}(A) = 1_{A}$,
the sums being bounded 
for $n^{th}$-order deformations, and extending to infinity for 
formal deformations.  In this case, the trivial deformation has
$\mu^{(j)} \equiv 0$ and $\iota^{(j)} \equiv 0$ for $j > 0$.

The definition then translates into equational conditions on the
$\mu^{(i)}$'s and $\iota^{(j)}$'s:

\begin{proposition}
The coefficients $\mu^{(i)}(f,g)$ and $\iota^{(j)}(A)$ 
in a 1-parameter deformation
of a category $\cal C$ satisfy

\begin{eqnarray}
\mu^{(0)}(f,g) & = & fg \label{eq:catdef1} \\
\iota^{(0)}(A) & = & 1_{A} \label{eq:catdef1i}\\
\sum_{i=0}^n \mu^{(i)}(\iota^{(n-i)}(s(f)),f) & = & 0 \,\, 
\mbox{\rm for $n >0$}
\label{eq:catdef2} \\
\sum_{i=0}^n \mu^{(i)}(f,\iota^{(n-i)}(t(f)) & = & 0 \,\, 
\mbox{\rm for $n >0$}
\label{eq:catdef3} \\
\sum_{i=0}^n \mu^{(n-i)}(\mu^{(i)}(a,b),c) & = & 
    \sum_{i=0}^n \mu^{(n-i)}(a,\mu^{(i)}(b,c)), \label{eq:catdef4}
\end{eqnarray}

\noindent and a family of coefficents $\{\mu^{(i)} | i=1,\ldots \}$ defines a
1-parameter deformation if and only if it satisfies these equations. 
\end{proposition}

\noindent{\sc Proof.} \ref{eq:catdef1} and \ref{eq:catdef1i} are
the requirement that 
the the deformation reduce to the identity modulo $\frak m = \langle
\epsilon
\rangle$. By trivial calculations \ref{eq:catdef2} and
\ref{eq:catdef3} are 
seen to be the
preservation of identities and \ref{eq:catdef4} the preservation 
of associativity. $\blacksquare$
\smallskip

Here, of course, there is an upper bound in the indices of the
$\mu^{(i)}$'s 
in the case of an $n^{th}$-order deformation, and no bound in the case
of a formal deformation.  We will not bother to note this again in the
discussion below of deformations of functors and natural
transformations.

Equivalences between 1-parameter deformations, similarly can be
characterized
in terms of coefficents.

In particular, we have

\begin{proposition}\label{catequiv} Given two 1-parameter 
    deformations of 
$\cal C$, $\hat{\cal C}_1$ and
$\hat{\cal C}_2$, with compositions given by 
$f\star_1 g = \sum \mu^{(i)}_1(f,g)\epsilon^{i}$ and 
$f\star_2 g = \sum
\mu^{(i)}_1(f,g)\epsilon^{i}$
respectively, and identity maps given by $\hat{1}_{1 A} =
\sum \iota^{(j)}_{1}(A)\epsilon^{j}$ and $\hat{1}_{2 A} = 
\sum \iota^{(j)}_{2}(A)\epsilon^{j}$
an equivalence $\Phi$ is given by 
a functor given on objects by the
identity, and on arrows by

\[ \Phi(f) = \sum \Phi^{(i)}(f)\epsilon^i, \] 

\noindent where 

\begin{eqnarray}
\Phi^{(0)}(f) & = & f \label{eq:catequiv1} \\
\sum_{j=0}^i \Phi^{(j)}(\iota^{(i-j)}_{1}(A)) & = & \iota^{(i)}_{2}(A) 
\;\; \mbox{\rm for $i>0$ and all $A \in Ob({\cal C})$} 
    \label{eq:catequiv2} \\
\sum_{j=0}^i \Phi^{(j)}(\mu_1^{(i-j)}(f,g)) & = & 
	\sum_{k+l+m = i} \mu_{2}^{(k)}(\Phi^{(l)}(f),\Phi^{(m)}(g)) 
        \label{eq:catequiv3}
\end{eqnarray}

\noindent Moreover, a family of assignments of parallel arrows  
to arrows in $\cal C$ defines an equivalence of 1-parameter
deformations 
if and only if it satisfies equations
\ref{eq:catequiv1}-\ref{eq:catequiv3}.
\end{proposition}

\noindent{\sc Proof.}  Equation \ref{eq:catequiv1} is the requirement
that $\Phi$ reduce to the identity functor modulo $\frak m$, equations
\ref{eq:catequiv2} and \ref{eq:catequiv3} are the preservation of
identity
arrows and composition respectively. Preservation of sources and
targets
is trivial, and invertibility follows immediately from the reduction
to the identity modulo $\frak m$. $\blacksquare$
\smallskip

Classical discussions of deformations of unital associative algebras 
sometimes omit or gloss over the question of deforming the identity 
element.  There is good reason for this, which holds with equal force 
in the many objects setting:

\begin{theorem}  \label{unitmightaswellbetrivial}
    If $\hat{\cal C}$ is a 1-parameter deformation of a 
    $k$-linear category $\cal C$, with deformed composition $f \star 
    g = \sum \mu^{(i)}(f,g)\epsilon^{i}$ and deformed identities $\hat{1}_{A} = 
    \sum \iota^{(j)}(A)\epsilon^{j}$, there exists an equivalent 
    deformation $\tilde{\cal C}$ with deformed composition $f \ast g 
    = \sum \tilde{\mu}^{(i)}(f,g)\epsilon^{i}$ and identity arrows equal to 
    the undeformed identities $1_{A}$.
    \end{theorem}
    
\noindent{\sc Proof.} We begin by constructing the functor from 
$\hat{\cal C}$ to $\tilde{\cal C}$:  map each 
object to itself, and each map $g$ to $1_{s(g)}\star g$.  Plainly
$\hat{1}_{A}$ is mapped to $1_{A}$.  

Now observe that in $\hat{\cal C}$, $1_{A}$ is an automorphism of 
$A$: 
the coefficients
$\kappa^{(i)}$ of its inverse $1_{A}^{-1}$
can be found inductively from the conditions $\kappa^{(0)} = 1_{A}$
and $\iota^{(n)} = \sum_{i=0}^{n} \mu^{(i)}(\kappa^{n-i},1_{A})$.

We can now define the composition on $\tilde{\cal C}$.  For 
$f:A\rightarrow B$ and $g:B\rightarrow C$, let $f \ast g 
= f \star 1_{B}^{-1} \star g$.  It is immediate by 
construction that $1_{A}$ is an identity for $\ast$.  Associativity of
$\ast$ follows from the associativity of $\star$, as does
the fact that the assignment of maps at the beginning of the proof
is a functor, since $1_{A}\star (f\star g) = (1_{A}\star f)\star
1_{B}^{-1} \star (1_{B} \star g)$. 

It is easy to see that reducing modulo $\langle \epsilon \rangle$
gives the identity functor of $\cal C$, and that   
the inverse functor is given on arrows by $f \mapsto 1_{A}^{-1}
\star f$. 

Likewise it is easy to see that $\tilde{\mu}^{(n)} = \sum_{i+j+k=n}
\mu^{(k)}(\mu^{j}(f,\kappa^{i}), g)$.
$\blacksquare$ 
\smallskip

For functors the natural notion of deformation is given by

\begin{definition}  A {\em deformation $\hat{F}$ of a $k$-linear functor
$F:{\cal C}\rightarrow {\cal D}$} is a triple $(\hat{\cal C},
\hat{\cal D}, \hat{F})$, where $\hat{\cal C}$ (resp. $\hat{\cal D}$)
is a deformation of $\cal C$ (resp. $\cal D$)
over the same local ring $R$, and $\hat{F}$ is a
functor enriched in the category of $R$-modules or $\frak m$-adically 
complete $R$-modules, as appropriate, from
$\hat{\cal C}$ to $\hat{\cal D}$ that reduces modulo $\frak m$
to $F$.
\end{definition}

As for categories, in the case of 1-parameter deformations, the
definition of
a deformation a $k$-linear functor is equivalent 
to a family of equational conditions on the coefficients of powers of 
$\epsilon$:

\begin{proposition}
If  $(\hat{\cal C}, \hat{\cal D}, \hat{F})$ is a 1-parameter
deformation of a 
functor $F:{\cal C}\rightarrow {\cal D}$ with the composition in 
$\hat{\cal C}$ (resp. $\hat{\cal D}$) given by 
$f \star g = \sum \mu^{(i})(f,g)\epsilon^i$ (resp. 
$f \ast g = \sum \nu^{(i)}(f,g)\epsilon^i$) and identities given
by $\hat{1}_{A} = \sum \iota^{(i)}(A)\epsilon^{i}$ (resp.
$\hat{1}_{A} = \sum \lambda^{(i)}(A)\epsilon^{i}$) and $\hat{F}$ given on
arrows 
by $\hat{F}(f) = \sum F^{(i)}(f)\epsilon^i$, then the $\mu^{(i)}$'s
and 
$\nu^{(i)}$'s satisfy equations \ref{eq:catdef1}-\ref{eq:catdef4},
and moreover

\begin{eqnarray}
F^{(0)} & = & F \label{eq:functdef1}  \\
\sum_{i+j=n}F^{(i)}(\iota^{(j)}(A)) & = & \lambda^{(n)}(F(A))
\,\, \makebox{\rm for $n>0$ and all $x \in
Ob({\cal C})$}
\label{eq:functdef2} \\
\sum_{i+j=n, \, i,j\geq 0} F^{(i)}(\mu^{(j)}(f,g)) & = & 
\sum_{k+l+m=n \, k,l,m\geq 0} \nu^{(k)}(F^{(l)}(f), F^{(m)}(g)),
\label{eq:functdef3}
\end{eqnarray}

\noindent and the families of coefficents define a 1-parameter
deformation
if and only if they satisfy these equations.

\end{proposition}

\noindent{\sc Proof.} \ref{eq:functdef1} is the requirement that 
$\hat{F}$ reduce modulo ${\frak m} = \langle \epsilon \rangle$ to
$F$.  
In the presence of \ref{eq:functdef1}
\ref{eq:functdef2} is equivalent to the preservation of identity
arrows by
$\hat{F}$, while \ref{eq:functdef3} is equivalent to the preservation
of
(the deformed) composition. $\blacksquare$
\smallskip 

Observe that an equivalence of deformations of a $k$-linear category 
$\cal C$ is simply a 
deformation of the identity functor that has the two deformations of
$\cal C$ as source and target.

There are evident notions of equivalence between deformations of
$k$-linear functors:

\begin{definition} Two deformations  $(\hat{\cal C},
\hat{\cal D}, \hat{F})$ and  $(\tilde{\cal C},
\tilde{\cal D}, \tilde{F})$ are {\em strongly equivalent}
if $\hat{\cal C} = \tilde{\cal C}$, $\hat{\cal D} = \tilde{\cal D}$, 
and there exists a natural isomorphism 
$\phi: \hat{F}\Rightarrow \tilde{F}$ which reduces to $Id_F$ modulo
$\frak m$.

Two deformations   $(\hat{\cal C},
\hat{\cal D}, \hat{F})$ and  $(\tilde{\cal C},
\tilde{\cal D}, \tilde{F})$ are {\em weakly equivalent} if there
exist equivalences of deformations of categories,
$\Gamma:\hat{\cal C}\rightarrow \tilde{\cal C}$ and 
$\Delta:\hat{\cal D}\rightarrow
\tilde{\cal D}$, and a natural 
isomorphism $\phi:\hat{F}\Delta \Rightarrow \Gamma \tilde{F}$,
which reduces modulo $\frak m$ to $Id_F$.
\end{definition}

For 1-parameter deformations, each type of equivalence can be 
characterized by equations 
on coefficients of powers of $\epsilon$.  We give the more general
case
of weak equivalence, as strong equivalence is simply specialization to
the case where $\Gamma$ and $\Delta$ are identity functors:

\begin{proposition} \label{functwe}
If $(\Gamma, \Delta, \phi)$ is a (weak) equivalence between 
1-parameter deformations $(\hat{\cal C}_1,\hat{\cal D}_1, \hat{F}_1)$
and  
$(\hat{\cal C}_2,\hat{\cal D}_2, \hat{F}_2)$ of 
$F:{\cal C}\rightarrow {\cal D}$, with 
$\Gamma(f) = \sum \Gamma^{(i)}(f)\epsilon^i$, 
$\Delta(f) = \sum \Delta^{(i)}(f)\epsilon^i$, 
and $\phi_x = \sum \phi^{(i)}_x\epsilon^i$, then the $\Gamma^{(i)}$'s 
(resp. $\Delta^{(i)}$'s) and the coefficients $\mu^{(i)}_n$ (resp. 
$\nu^{(i)}_n$) defining the composition
on $\hat{\cal C}_n$ (resp. $\hat{\cal D}_n$) $n = 1,2$ satisfy 
equations \ref{eq:catequiv1}-\ref{eq:catequiv3} {\em mutatis
mutandis},
and

\begin{eqnarray}
\phi^{(0)}_x & = & 1_x \,\, \mbox{\rm for all $x \in Ob({\cal C})$} 
      \label{eq:functequiv1} \\
\sum_{i+j+k=n} \Delta^{(i)}(F^{(j)}_1(f))\phi^{(k)}_y & = &
      \sum_{p+q+r=n} \phi^{(p)}_x F^{(q)}_2(\Gamma^{(r)}(f))
      \label{eq:functequiv2} \\
     & & \mbox{\rm for all $n \geq 0$ and all 
        $f:x\rightarrow y \in Arr({\cal C})$}. \nonumber
\end{eqnarray}

\noindent Moreover a family of coefficients defines a weak equivalence of 
1-parameter deformations if and only if it satisfies the given
conditions. 

\end{proposition} 

\noindent{\sc Proof.} The conditions not involving $\phi$ must hold
by Proposition
\ref{catequiv}.  Equation \ref{eq:functequiv1} is the requirement that
$\phi$ reduce to the identity natural transformation modulo $\frak
m$.  
Equation \ref{eq:functequiv2} is naturality. $\blacksquare$ \smallskip

Theorem \ref{unitmightaswellbetrivial} can be extended to functors:

\begin{theorem} \label{trivialunitsforfunctors} 
    If $F:{\cal A}\rightarrow {\cal B}$ is a functor, and
   $(\hat{\cal A},\hat{\cal B},\hat{F})$ is a deformation, then there 
   is a weakly equivalent deformation in which the identity arrows of
   ${\cal A}$ and ${\cal B}$ are undeformed.
   \end{theorem}
   
\noindent{\sc Proof.}  Construct equivalent deformations $\tilde{\cal 
A}$ and $\tilde{\cal B}$ of the 
categories as in Theorem \ref{unitmightaswellbetrivial}.  Let
$\Phi_{\cal X}:\hat{\cal X}\rightarrow \tilde{\cal X}$, 
$\Psi_{\cal X}:\tilde{\cal X}\rightarrow \hat{\cal X}$, 
$\sigma_{\cal X}:\Phi_{\cal X}\Psi_{\cal X} \Rightarrow Id_{\hat{\cal C}}$,
and $\tau_{\cal X}:\Psi_{\cal X}\Phi_{\cal X} \Rightarrow 
Id_{\tilde{\cal C}}$ define the equivalence of categories for ${\cal 
X} = {\cal A}, {\cal B}$.

Then $\tilde{F} := \Psi_{A}\hat{F}\Phi_{B}$ gives the deformed functor from
$\tilde{\cal A}$ to $\tilde{\cal B}$, and $\Phi_{\cal 
B}(F(\sigma_{\cal A}))$ gives the necessary natural isomorphism.
(Recall our convention that compositions are in diagrammatic order 
unless parentheses indicate `evaluation at'.)
$\blacksquare$\smallskip

Finally for natural transformations, we make

\begin{definition}
A {\em deformation $\hat{\sigma}$ of a natural 
transformation $\sigma:F\Rightarrow G$} 
(for functors $F,G:{\cal C}\rightarrow {\cal D}$ is a five-tuple
$(\hat{\cal C}, \hat{\cal D}, \hat{F}, \hat{G}, \hat{\sigma})$, where 
$\hat{\cal C}$ (resp. $\hat{\cal D}, \hat{F}, \hat{G}$) is a
deformation of 
$\cal C$ (resp. ${\cal D}, F, G$), and $\hat{\sigma}$ is a natural 
transformation from $\hat{F}$ to $\hat{G}$, which reduces modulo
$\frak m$
to $\sigma$.
\end{definition}

For deformations of natural transformations, there are three obvious
notions 
of equivalence:

\begin{definition}
Two deformations $(\hat{\cal C}, \hat{\cal D}, \hat{F}, \hat{G},
\hat{\sigma})$
and $(\tilde{\cal C}, \tilde{\cal D}, \tilde{F}, \tilde{G},
\tilde{\sigma})$ 
of a natural transformation $\sigma$ are {\em strongly equivalent} if 
$(\hat{\cal C}, \hat{\cal D}, \hat{F}, \hat{G}) = 
(\tilde{\cal C}, \tilde{\cal D}, \tilde{F}, \tilde{G})$.

Two deformations $(\hat{\cal C}, \hat{\cal D}, \hat{F}, \hat{G},
\hat{\sigma})$
and $(\tilde{\cal C}, \tilde{\cal D}, \tilde{F}, \tilde{G},
\tilde{\sigma})$ 
are
{\em equivalent} if $\hat{\cal C} = \tilde{\cal C}$,  $\hat{\cal D} = 
\tilde{\cal D}$, and $\hat{F}$ (resp. $\hat{G}$) is strongly
equivalent to 
$\tilde{F}$ (resp. $\tilde{G}$) by a strong equivalence of
deformations of
functors $\phi$ (resp. $\gamma$), and, moreover, $\hat{\sigma}\gamma =
\phi\tilde{\sigma})$.

Finally, two deformations 
$(\hat{\cal C}, \hat{\cal D}, \hat{F}, \hat{G}, \hat{\sigma})$
and $(\tilde{\cal C}, \tilde{\cal D}, \tilde{F}, \tilde{G},
\hat{\sigma})$ are
{\em weakly equivalent} if there exist equivalences of deformations 
of categories, $\Gamma:\hat{C}\rightarrow \tilde{C}$ and
$\Delta:\hat{D}\rightarrow \tilde{D}$, and a natural 
isomorphisms $\phi:\hat{F}\Delta \Rightarrow \Gamma \tilde{F}$ (resp.
$\psi:\hat{G}\Delta \Rightarrow \Gamma \tilde{G}$),
which reduces modulo $\frak m$ to $Id_F$ (resp. $Id_G$), and, moreover
$\hat{\sigma}\psi =
\phi\tilde{\sigma})$.
\end{definition}

Here again, in the case of 1-parameter deformations, the definition
is equivalent to a family of equational conditions on the coefficients
of $\epsilon$:

\begin{proposition}
If $(\hat{\cal C}, \hat{\cal D}, \hat{F}, \hat{G}, \hat{\sigma})$ is
a 1-parameter deformation of a natural transformation
$\sigma:F\Rightarrow G$
with the composition on 
$\hat{\cal C}$ (resp. $\hat{\cal D}$) given by 
$f \star g = \sum \mu^{(i})(f,g)\epsilon^i$ (resp. 
$f \star g = \sum \nu^{(i)}(f,g)\epsilon^i$), $\hat{F}$ (resp.
$\hat{G}$) given on arrows 
by $\hat{F}(f) = \sum F^{(i)}(f)\epsilon^i$ (resp. 
$\hat{G}(f) = \sum G^{(i)}(f)\epsilon^i$), and
$\hat{\sigma}_x = \sum \sigma^{(i)}_x \epsilon^i$ then the
$\mu^{(i)}$'s and 
$\nu^{(i)}$'s satisfy equations \ref{eq:catdef1}-\ref{eq:catdef4},
they and the $F^{(i)}$'s and $G^{(i)}$'s satisfy 
\ref{eq:functdef1}-\ref{eq:functdef3}, and moreover

\begin{eqnarray}
\sum_{i+j+k=n \, i,j,j\geq 0} \nu^{(i)}(F^{(j)}(f),\sigma^{(k)}_y)
& = & \sum_{p+q+r=n \, p,q,r\geq 0}
\nu^{(p)}(\sigma^{(q)}_x,G^{(r)}(f))
\label{eq:nattransdef}
\end{eqnarray}

\noindent for all $f:x\rightarrow y$ in $\cal C$.

Moreover families of coefficients define a 1-parameter deformation
if and only if they satisfy these conditions.
\end{proposition}

\noindent{\sc Proof.} The conditions not involving the
$\sigma^{(i)}$'s
must hold by earlier propositions. Equation \ref{eq:nattransdef} is 
simply the naturality of $\hat{\sigma}$. $\blacksquare$
\smallskip

And finally, equivalences between 
1-parameter deformations of natural transformations
can be reduced to equations on the coefficients:

\begin{proposition}
If $(\Gamma, \Delta, \phi, \psi)$ is a weak equivalence between
two deformations of a natural transformation $\sigma:F\Rightarrow G$
between functors $F,G:{\cal C}\rightarrow {\cal D}$, 
$(\hat{\cal C}_1, \hat{\cal D}_1, \hat{F}_1, \hat{G}_1,
\hat{\sigma}_1)$
and $(\hat{\cal C}_2, \hat{\cal D}_2, \hat{F}_2, \hat{G}_2,
\hat{\sigma}_2)$, with 
$\Gamma(f) = \sum \Gamma^{(i)}(f)\epsilon^i$, 
$\Delta(f) = \sum \Delta^{(i)}(f)\epsilon^i$, 
$\phi_x = \sum \phi^{(i)}_x\epsilon^i$, $\psi_x = \sum
\psi^{(i)}_x\epsilon^i$,
then
the coefficents defining the compositions on the $\hat{\cal C}_n$'s 
and $\hat{\cal D}_n$'s ($n = 1,2$) together with the $\Gamma^{(i)}$'s
and
the $\Delta{(i)}$'s, the $F^{(i)}_n$'s ($n = 1,2$) (resp.
$G^{(i)}_n$'s ($n= 1,2$)) and the $\phi^{(i)}$'s (resp. the
$\psi^{(i)}$'s) 
satisfy the conditions of Proposition \ref{functwe} {\em mutatis 
mutandis}, and moreover

\begin{eqnarray} \label{eq:nattranswe} 
\sum_{i+j=n} \phi^{(i)}_x\sigma^{(j)}_{2 \; x} =
\sum_{p+q=n} \sigma^{(p)}_{1 \; x}\psi^{(q)}_x
\end{eqnarray}

\noindent for all $n$ and $x \in Ob({\cal C})$.

Moreover, families of coefficients define a weak equivalence of
deformations
of natural transformations if and only if they satisfy these
conditins.
\end{proposition}

\noindent{\sc Proof.}
All but the last condition must hold by Proposition \ref{functwe} and
the definition of weak equivalence.  The last is simply the condition
that
$\sigma_1\psi = \phi\sigma_2$.
$\blacksquare$
\smallskip

Again, up to weak equivalence, we may assume that identity maps have 
been left undeformed, as the following follows by the construction of 
Theorem \ref{trivialunitsforfunctors}, and a little 2-categorical 
diagram chase around a `square pillow' with two squares and
two bigons as faces, the equivalences from
Theorem \ref{trivialunitsforfunctors} as the seams between the square 
faces, the square faces being the natural transformation in the weak 
equivalences between the deformations of $F$ and $G$, and one bigon 
being the original deformation of $\sigma$:

\begin{theorem}
    If $(\hat{\cal C}, \hat{\cal D}, \hat{F}, \hat{G},\hat{\sigma})$ 
    is a deformation of a natural transformation $\sigma:F 
    \Rightarrow G$ for $F,G:{\cal C}\rightarrow {\cal D}$, then there 
    exists a weakly equivalent deformation of $\sigma$ in which the
    identity maps of both $\cal C$ and $\cal D$ are undeformed.
    \end{theorem}

Having arrived at this point, it is clear what a deformation of a
3-cell is:  it is simply a deformation of the equal bounding 2-cells.

In the case where the diagram is simply a `bigonal pillow', that is a
3-computad with a single 3-cell with single 2-cells as source and
target, there is nothing more to be said.
We will see, however, once we consider pasting of deformations, and
deformations of pasting diagrams, that in general the matter is not
trivial.

\section{The Hochshild cohomology of $k$-linear categories,\\
functors and natural transformations: definitions and elementary 
results \label{Hochschild}}

As was the case in \cite{GSdiag}, the cohomology appropriate to
objects turns out 
to be a special case of that appropriate to arrows, though with
groups in 
different cohomological dimensions playing corresponding role.  Here,
however,
we can apply this observation twice. 

It turns out that the obvious Hochschild complex for a natural 
transformation depends only on its source and target, we thus begin
with

\begin{definition}
The {\em Hochschild complex of a parallel pair of functors, }

\[F,G:{\cal A}\rightarrow
{\cal B}\]

\noindent has cochain groups given by

\[ C^n(F,G) := \prod_{x_0,\ldots ,x_n\in Ob({\cal A})}
  {\rm Hom}_k ({\cal A}(x_0,x_1)\otimes \ldots \otimes {\cal
A}(x_{n-1},x_n),
{\cal B}(F(x_0), G(x_n))), \]

\noindent and for $n = 0$, noting that $k$ is the empty tensor 
product,

\[C^{0}(F,G) := \prod_{x_{0}\in Ob({\cal A})} {\rm Hom}_{k}(k, {\cal 
B}(F(x_{0}), G(x_{0}))) \]

\noindent with coboundary given by

\begin{eqnarray*}
  \lefteqn{\delta \psi (f_0\otimes \ldots \otimes f_n) :=}\\
  & & F(f_0)\phi(f_1\otimes \ldots \otimes f_n) +\\ 
  & & \sum_{i=1}^n (-1)^i 
\phi(f_0\otimes \ldots \otimes f_{i-1}f_i \otimes \ldots \otimes
f_n) +\\ 
  & & (-1)^{n+1} \phi(f_0\otimes \ldots \otimes f_{n-1})G(f_n) 
\end{eqnarray*}
 
\end{definition}

$\delta^2 = 0$ by the usual calculation.

The impression that natural transformations themselves are forgotten 
in this definition is deceptive.  In fact, in this context naturality 
itself turns out to be a cohomological condition:

\begin{proposition}
    A natural transformation from $F$ to $G$ is a 0-cocycle in 
    $C^\bullet (F,G)$, or equivalently a 0-dimensional cohomology 
    class.
    \end{proposition}
    
\noindent {\sc Proof.} Observe that 

\[ C^{0}(F,G) = \prod_{x \in 
Ob({\cal A})} {\rm Hom}_{k}(k, {\cal B}(F(x),G(x))). \]  

\noindent (In the 
sources, one has an empty tensor products of hom-spaces in $\cal A$,
since there is no ``next'' element object, while in the target one 
has the hom-space from the image of the first object under $F$ to the 
last, here the same object, under $G$.)

Thus a 0-cochain is an assignment of an arrow 
$\phi_{x}:F(x)\rightarrow G(x)$ 
for each object $x$ of $\cal A$.

The cocycle condition then is

\[ 0 = \delta(\phi)(f:x_{0}\rightarrow x_{1}) = 
F(f)\phi_{x_{1}} - \phi_{x_{0}}G(f), \]

\noindent that is, the naturality of $\phi$.  Since
$C^{-1}(F,G) = 0$, the 0-cohomology classes and 0-cocycles can
be identified. $\blacksquare$

We will also be interested in the subcomplex of normalized Hochschild
cochains--those which vanish whenever one of the arguments is an
identity arrow.  We denote this subcomplex by $\bar{C}^\bullet(F,G)$.

Theorem \ref{unitmightaswellbetrivial} has a cohomological analog:

\begin{theorem} The normalized Hochschild comples $\bar{\cal 
    C}^{\bullet}(F,G)$ is a chain deformation retract of the 
    Hochschild complex ${\cal C}^{\bullet}(F,G)$.
    \end{theorem}
    
\noindent {\sc Proof.} The result follows from the same trick used in 
\cite{Laz}. 

Call a cochain $\phi$ $i$-normalized if $\phi(f_{0},\ldots ,f_{n-1})$ 
is zero whenever $f_{j}$ is an identity arrow for any $j \leq i$.

The $i$-normalized cochains from a subcomplex ${\cal 
C}^{\bullet}_{i}(F,G)$ of ${\cal C}^{\bullet}(F,G)$, and satisfy

\[{\cal C}^{\bullet}_{i+1}(F,G) \subset {\cal 
C}^{\bullet}_{i}(F,G), \]

\noindent and 

\[ \cap_{i=0}^{\infty} {\cal 
C}^{\bullet}_{i}(F,G) = \bar{\cal C}^{\bullet}(F,G). \]

For $k \geq 0$, define maps $s^{k}:{\cal 
C}^{n}(F,G)\rightarrow {\cal C}^{n-1}(F,G)$ by

\[ s^{k}(\phi)(f_{1},\ldots f_{n-1}) = \left\{ \begin{array}{ll}
0 & \mbox{if $k > n$} \\
\phi(f_{1},\ldots, f_{k}, 1_{t(f_{k})}, f_{k+1},\ldots ,f_{n-1}) &
\mbox{if $k \leq n$} \end{array} \right. \]

Let $h^{k}(\phi) := \phi - \delta(s^{k}(\phi) - s^{k}(d(\phi))$.

Then $h^{0}$ and the inclusion $i_{0}:{\cal 
C}^{\bullet}_{0}(F,G)\rightarrow {\cal C}^{\bullet}(F,G)$ form
a chain deformation retraction of the whole complex onto the
subcomplex, with $s^{0}$ as the homotopy from $i_{0}h^{0}$ to 
the identity of ${\cal C}^{\bullet}(F,G)$.
And, $h^{k}$ and the inclusion of
$i_{k}:{\cal C}^{\bullet}_{k-1}(F,G)\rightarrow 
{\cal C}^{\bullet}_{k}(F,G)$
form a chain deformation retraction, with $s^{k}$ as the homotopy
from $i_{k}h^{k}$ to the identity of ${\cal C}^{\bullet}_{k-1}(F,G)$.

Note that $h^{k}$ is the identity map on ${\cal C}^{n}(F,G)$ for
$n < k$.  Thus, the formal infinite composition 
$\bar{h} = h^{0}h^{1}h^{2}\ldots$ is finite in each dimension, and gives
a chain deformation retraction of ${\cal C}^{\bullet}(F,G)$ onto
$\bar{\cal C}^{\bullet}(F,G)$, with $s^{0} + h^{0}s^{1}
+ h^{0}h^{1}s^{2} + \ldots$ as the homotopy from $\bar{i}\bar{h}$
to the identity of ${\cal C}^{\bullet}(F,G)$, where
$\bar{i}$ is the inclusion of the subcomplex of normalized cochains.
$\blacksquare$

We also make

\begin{definition}
The {\em Hochschild complex of a functor $F$} is 

\[ (C^\bullet (F),\delta) := 
(C^\bullet(F,F),\delta). \]  

The {\em Hochschild complex of a $k$-linear
category $\cal C$} is 

\[ (C^\bullet({\cal C}), \delta) := (C^\bullet(Id_{\cal C},
Id_{\cal C}), \delta). \]
\end{definition}

The subcomplexes of normalized cochains 
$\bar{C}^\bullet(F)$ and
$\bar{C}^\bullet({\cal C})$ are defined in the obvious way.

Observe that this definition agrees with that given in \cite{Anel}
and in the 
special case of algebras with that given in \cite{G1, G2}.

Not surprisingly, the Hochschild complexes admit a rich algebraic
structure generalizing
that discovered by Gerstenhaber \cite{G1, G2} on the Hochschild
complex of 
an associative algebra (cf. also \cite{GV}).

First, given three parallel functors $F,G,H:{\cal A}\rightarrow {\cal
B}$,
there is a cup-product, or more properly a cup-product-like
2-dimensional
(in the sense of bicategory theory) composition,
$\cup:C^n(F,G)\otimes C^m(G,H)\rightarrow C^{n+m}(F,H)$ given by

\[ \phi \cup \psi (f_1\otimes \ldots \otimes f_{n+m}) :=
(-1)^{nm} \phi(f_1\otimes \ldots \otimes f_n)\psi(f_{n+1}\otimes
\ldots
\otimes f_{n+m}) \]

Second, given functors $G,H:{\cal B}\rightarrow {\cal C}$ and
$F_0,\ldots, F_n:{\cal A}\rightarrow {\cal B}$, 
there is a brace-like 1-dimensional composition 
(cf. \cite{GV})
\medskip

\centerline{\parbox{9cm}{\raggedright{
$ -\{-,\ldots,-\}:C^K(G,H)\otimes C^{k_1}(F_0,F_1) \otimes \ldots
\otimes C^{k_n}(F_{n-1},F_n)\rightarrow$
\vspace*{3pt}
\hspace*{3cm}$ C^{K + k_1 + \ldots k_n - n}(G(F_0),H(F_n))$ }}}
\medskip
 
\noindent given by

\begin{eqnarray*} 
\lefteqn{ \phi\{\psi_1,\ldots , \psi_n\}(f_1,\ldots ,f_N) :=}\\
 & & \sum (-1)^\epsilon
\phi(F_0(f_1), \ldots, \psi_1(f_{l_1+1}, \ldots ,f_{l_1+k_1}), 
F_1(f_{l_1+k_1+1}),\\
 & & \ldots \psi_2(f_{l_2+1}, \ldots ,f_{l_2+k_2}), \ldots,
\psi_n(f_{l_n+1}, \ldots , f_{l_n+k_n}), F_n(f_{l_n+k_n+1}), \ldots, F_n(f_N))
\end{eqnarray*}

\noindent where $N = K + k_1 + \ldots k_n - n$ and in each term
$\epsilon =
\sum_{i=1}^n (k_i - 1)l_i$, where $l_i$ is the total number of inputs 
occurring before $\psi_i$, and the outer sum ranges over all
insertions of
the $\psi_i$'s, in the given order, with any number, including zero,
of the
arguments preceeding $\psi_1$ (resp. between $\psi_i$ and
$\psi_{i+1}$, 
following $\psi_n$) with $F_0$ (resp. $F_i$, $F_n$) applied, and a
total of
$N$ arguments, including both those inside and outside of the
$\psi_i$'s.

In what follows, we will have call to consider a special instance of 
the brace-like 1-composition:  given a pair of parallel functors 
$F,G:{\cal A}\rightarrow {\cal B}$, there is a map (of graded vector
spaces, though not of cochain complexes)

\[ -\{-\}: C^{\bullet+1}({\cal B})\otimes C^{0}(F,G)\rightarrow 
C^{\bullet}(F,G) \]

The brace-like 1-composition and cup-product-like 2-composition
satisfy 
the 
identities given in Gerstenhaber and Voronov \cite{GV} whenever
sources and
targets are in agreement so that both sides of the equation are
defined.  
Also note that if
$n > K$ $\phi\{\psi_1,\ldots \psi_n\} = 0$, since the outer 
sum in the definition is empty.

Of course, in the many objects setting, it is not in general
possible
to reverse the order of $x\{y\}$, except in special cases, nor to add
the two
orderings, without yet more conditions.  The most general case in
which 
the usual construction gives rise to a differential graded Lie
algebra 
structure on $C^\bullet (F,G)$ is that where both $F$ and $G$ are
idempotent
endofunctors on some category. This includes, of course, the special case
$C^\bullet(Id_{\cal C},Id_{\cal C}) = C^\bullet({\cal C})$.

The circumstance for typical pairs of parallel functors, however, 
suggests that the 
view that deformation theory is governed by a differential graded
Lie algebra ought to be generalized to regard the brace algebra
structure
as more fundamental.  In view of this observation, we will consider 
analogs of the Maurer-Cartan equation for functors and natural 
transformations expressed in terms of
the brace-like composition rather than seeking a formulation in terms
of
a dg-Lie algebra or $L^\infty$-algebra structure.

There are also a number of rather obvious cochain maps between these 
complexes induced by the usual 2-categorical operations:

\begin{proposition} \label{cochain1comp}
If $F,G:{\cal C}\rightarrow {\cal D}$, and $H:{\cal D}\rightarrow
{\cal E}$
are functors, then there is a cochain map 
$H_*(-):C^\bullet(F,G)\rightarrow
C^\bullet(H(F),H(G))$
given by 

\[ H_*(\phi)(f_1,\ldots ,f_n) := H(\phi(f_1,\ldots ,f_n)). \]

Similarly if $F:{\cal C}\rightarrow {\cal D}$ and 
$G,H:{\cal D}\rightarrow {\cal E}$ are functors, there is a cochain
map
$F^*(-) = -(F^\bullet):C^\bullet(G,H)\rightarrow
C^\bullet(G(F),H(F))$ given by

\[ F^*(\phi)(f_1,\ldots , f_n) = \phi(F^\bullet)(f_1,\ldots , f_n) := 
\phi(F(f_1),\ldots ,F(f_n)). \]
\end{proposition}  

\noindent{\sc Proof.} When it is remembered that all functors are
linear,
and that the analog of the left and right actions on the bimodule of
coefficients involve the source and target functors of the pair, both
statements follow by trivial calculations.
$\blacksquare$  

Both of these cochain maps can be expressed in terms of the 
brace-like 1-composition:  for any functor $K:{\cal X}\rightarrow
{\cal Y}$ the action of $K$ on the hom vectorspaces gives a family
of linear maps ${\cal X}(x,y)\rightarrow {\cal Y}(K(x),K(y))$, 
defining a $1$-cochain ${\Bbb K}\in C^{1}(K,K)$.  

We then have

\[ H_{*}(\phi) = {\Bbb H}\{\phi\} \]

\noindent and

\[ F^{*}(\phi) = \phi\{{\Bbb F},\ldots ,{\Bbb F}\}. \]

\smallskip
Notice in the special case where ${\cal C} = {\cal D}$ and 
$F = G = Id_{\cal C}$, $H_*$ defines a cochain map 
$H_*: C^\bullet({\cal C})\rightarrow C^\bullet(H)$, while in the case
where ${\cal D} = {\cal E}$ and $G = H = Id_{\cal D}$, $F^*$ defines
a cochain map $F^*: C^\bullet({\cal D})\rightarrow C^\bullet(F)$.

Of particular interest is the cochain map 

\[ F^*(p_2) - F_*(p_1): C^\bullet({\cal C})\oplus C^\bullet({\cal D})
\rightarrow C^\bullet(F) \]

\noindent the cone over which will occur in the classification of
deformations
of functors. 

Equally trivial is the proof of 

\begin{proposition} \label{cochain2comp}
If $\tau:F_1\Rightarrow F_2$ is a natural transformation, then
post- (resp. pre-) composition by $\tau$ induces a cochain map
$\tau^* :C^\bullet(F_2,G)\rightarrow C^\bullet(F_1,G)$ 
(resp. $\tau_*:C^\bullet(G,F_1)\rightarrow C^\bullet(G,F_2)$
for any
functor $G$. 
\end{proposition}

These cochain maps can be expressed in terms of the cup-product-like
2-composition:

\[ \tau^{*}(\phi) = \tau \cup \phi \]

\noindent and

\[ \tau_{*}(\phi) = \phi \cup \tau. \]

These cochain maps, however, are less important than one induced by
a natural transformation $\sigma:F\Rightarrow G$ 
from the cone over the
map 

\[ \left[ \begin{array}{cc}
                  F_* & -F^* \\               
                  G_* & -G^*  
               \end{array} \right]: 
             C^\bullet({\cal A})\oplus C^\bullet({\cal B})
               \rightarrow C^\bullet(F) \oplus C^\bullet(G). \]

\noindent to $C^\bullet(F,G)$. 	       

Let $C^\bullet({\cal A} \stackrel{F}{G} {\cal B})$ denote this cone,
so the cochain groups are

\[ C^\bullet({\cal A} \stackrel{F}{G} {\cal B}) :=
C^{\bullet+1}({\cal A})
\oplus C^{\bullet+1}({\cal B})\oplus C^\bullet(F)\oplus C^\bullet(G)
\]

with coboundary operators given by

\[ d_{{\cal A} \stackrel{F}{G} {\cal B}} = \left[
    \begin{array}{cccc}
	-d_{\cal A} & 0 & 0 & 0 \\
	0 & -d_{\cal B} & 0 & 0\\
	-F_{*} & F^{*} & d_{F} & 0\\
	-G_{*} & G^{*} & 0 & d_{G}
     \end{array} \right] \]	

We then have:

\begin{proposition} Let $\sigma:F\Rightarrow G$ be a natural
transformation, 
then

\[ \sigma \ddagger := 
\left[ \begin{array}{cccc} 0 & (-)\{\sigma \} & \sigma_*
&  \sigma^* \end{array} \right]: 
C^\bullet({\cal A} \stackrel{F}{G} {\cal B}) \rightarrow
C^\bullet(F,G) \]

\noindent is a cochain map.
\end{proposition}

\noindent{\sc Proof.}  Recall that a natual transformation is a 
0-cocycle in $C^{\bullet}(F,G)$, and thus we can use the brace-like
1-composition to define the second entry.

Using subscripts to distinguish between the various coboundaries, we 
then have

\[ d_{F,G}\left(\sigma \ddagger\left(\left[ \begin{array}{c}\phi \\ \psi \\
\upsilon \\ \omega \end{array}\right]\right)\right) = 
d_{F,G}(\psi\{\sigma\}) + d_{F,G}(\sigma_{*}(\upsilon)) -
d_{F,G}(\sigma^{*}(\omega)) \]

\noindent and

\begin{eqnarray*} 
    \sigma \ddagger \left(d_{{\cal A} \stackrel{F}{G} {\cal B}}\left(
\left[\begin{array}{c}\phi \\ \psi \\
\upsilon \\ \omega \end{array}\right]\right)\right)  & = & 
 -d_{\cal B}(\psi)\{\sigma\} - \sigma_{*}(F_{*}(\phi)) + 
\sigma_{*}(F^{*}(\psi))  \\ & & + \sigma_{*}(d_{F}(\upsilon) +  
\sigma^{*}(G_{*}(\phi)) - \sigma^{*}(G^{*}(\psi)) 
- \sigma^{*}(d_{G}(\omega), \end{eqnarray*}

\noindent which we must show are equal.  

The terms involving $\phi$
in the second expression cancel by the naturality of $\sigma$, while 
ther terms involving $\upsilon$ (resp. $\omega$) in the two expressions 
are equal since $\sigma_{*}$ (resp. $\sigma^{*}$) is a cochain map.
It thus remains only to show that the terms involving $\psi$ agree.

Expanding $d_{\cal B}(\psi)\{\sigma\}(f,g)$ for 
$X\stackrel{f}{\rightarrow}Y\stackrel{g}{\rightarrow}Z$ gives

\[ \sigma_{X}\psi(G(f),G(g) - \psi(\sigma_{X}G(f),G(g)) + 
\psi(\sigma_{X},G(fg)) - \psi(\sigma_{X},G(f))G(g) \]
\[ -F(f)\psi(\sigma_{Y},G(g)) + \psi(F(f)\sigma_{Y},G(g)) 
- \psi(F(f),\sigma_{Y}(G(g)) + \psi(F(f),\sigma_{Y})G(g)  \]
\[ + F(f)\psi(F(g),\sigma_{Z}) - \psi(F(fg),\sigma_{Z}) + 
\psi(F(f),F(g)\sigma_{Z}) - \psi(F(f),F(g))\sigma_{Z} \]

\noindent  The terms in which one of the arguments of $\psi$ consists
of an instance of $\sigma$ composed with another map (in either 
order) cancel in pairs by the naturality of $\sigma$.  The first and 
last term are $\sigma^{*}(G^{*}(\psi))(f,g)$ and 
 $-\sigma_{*}(F^{*}(\psi))(f,g)$  respectively.  The remaining terms are 
easily seen to be $-d_{F,G}(\psi\{\sigma\})(f,g)$, yielding the 
desired result. $\blacksquare$
\smallskip

Because of its importance to the deformation theory of natural
transformations,
we will denote the cone on $\sigma \ddagger$ by ${\frak
C}^\bullet(\sigma)$,
likewise we will denote the cone on $F^*(p_2) - F_*(p_1)$ by 
${\frak C}^\bullet(F)$, and for completeness, we will let
${\frak C}^\bullet({\cal C})$ be another notation for the Hocschild
complex of ${\cal C}$.

More explicitly, 

\[{\frak C}^\bullet(\sigma) = C^{\bullet+2}({\cal A})
\oplus C^{\bullet+2}({\cal B})\oplus C^{\bullet+1}(F)\oplus 
C^{\bullet+1}(G)\oplus C^{\bullet}(F,G) \]

\noindent with coboundary operators given by

\[ {\frak d}_{\sigma} = \left[
    \begin{array}{ccccc}
	d_{\cal A} & 0 & 0 & 0 & 0\\
	0 & d_{\cal B} & 0 & 0 & 0\\
	F_{*} & -F^{*} & -d_{F} & 0 & 0\\
	G_{*} & -G^{*} & 0 & -d_{G} & 0\\
	0 & -(-)\{\sigma \} & -\sigma_* &  -\sigma^* & d_{F,G}
     \end{array} \right] .\]
     
While

\[ {\frak C}^\bullet(F) = C^{\bullet+1}({\cal A})\oplus 
C^{\bullet+1}({\cal B})\oplus C^{\bullet}(F) \]

\noindent with coboundary operators given by

\[ {\frak d}_{F} = \left[
    \begin{array}{ccc}
	-d_{\cal A} & 0 & 0 \\
	0 & -d_{\cal B} & 0 \\
	F_{*} & -F^{*} & d_{F}
	\end{array} \right] .\]

We will refer to ${\frak C}^\bullet(\sigma)$ (resp.
${\frak C}^\bullet(F)$, ${\frak C}^\bullet({\cal C})$) as {\em the
deformation
complex of the natural transformation (resp. functor, category)}, and 
denote its cohomology groups by  ${\frak H}^\bullet(\sigma)$ (resp.
${\frak H}^\bullet(F)$, ${\frak H}^\bullet({\cal C})$).
 
It is the purpose of the next section to justify these names by
showing
that first order deformations are classified up to equivalence by the
expected cohomology group, and that the obstructions to extending a
deformation to higher order all 
lie in the next cohomological dimension.

Finally, we introduce a deformation complex for a 3-cell
$1_\sigma:\sigma \equiv\!\rangle \sigma$.  It might seem reasonable to have
this simply be ${\frak C}^\bullet(\sigma)$ again.  However, for
reasons which will become clear once begin considering pasting
diagrams in general, it will be better to use a weakly equivalent
complex:

First let $\hat{\frak C}^\bullet(\sigma)$ denote the cone
on 

\[ i_1(\sigma\ddagger) + i_2(\sigma\ddagger):
C^\bullet({\cal A} \stackrel{F}{G} {\cal B}) \rightarrow
C^\bullet(F,G)\oplus
C^\bullet(F,G) ,\]

Then ${\frak C}^\bullet(1_\sigma)$ is the cone on 

\begin{eqnarray*} \lefteqn{p_5 - p_6: \hat{\frak C}^\bullet(\sigma)
=}\\
& & C^{\bullet+2}({\cal A})
\oplus C^{\bullet+2}({\cal B})\oplus C^{\bullet+1}(G)\oplus
C^{\bullet+1}(F)
\oplus C^\bullet(F,G) \oplus C^\bullet(F,G)\\
& & \rightarrow C^\bullet(F,G). \end{eqnarray*}

Observe that this is plainly weakly equivalent to ${\frak
C}^\bullet(\sigma)$.

\section{First order deformations 
without pasting and cohomology}\label{nopasting}

Let us begin by considering the first order case of the equational
conditions
defining 1-parameter deformations.  For categories and functors, we
obtain
the obvious generalization of the results of Gerstenhaber and Schack
\cite{GSdiag} to the many-objects case:

For categories, equations \ref{eq:catdef1}-\ref{eq:catdef4} become

\begin{eqnarray}
\mu^{(0)}(f,g) & = & fg \label{eq:catdef1.1} \\
\mu^{(1)}(1_{s(f)},f) & = & 0 \label{eq:catdef2.1} \\
\mu^{(1)}(f,1_{t(f)}) & = & 0 \label{eq:catdef3.1} \\
\mu^{(1)}(a,b)c + \mu^{(1)}(ab,c) & = & 
    a\mu^{(1)}(b,c) + \mu^{(1)}(a,bc) \label{eq:catdef4.1}
\end{eqnarray}

Equation \ref{eq:catdef4.1}, as expected, says that $\mu^{(1)}$ is a
Hochschild
2-cocycle, while equations \ref{eq:catdef2.1} and \ref{eq:catdef3.1}
require it to be normalized in the obvious sense.

The equations defining an equivalence of such deformations,
\ref{eq:catequiv1}-\ref{eq:catequiv3} become

\begin{eqnarray}
\Phi^{(0)}(f) & = & f \label{eq:catequiv1.1} \\
\Phi^{(1)}(1_x) & = & 0 \;\; \mbox{\rm for all $x \in Ob({\cal C})$} 
    \label{eq:catequiv2.1} \\
\mu_1^{(1)}(f,g) + \Phi^{(1)}(fg) & = & 
\mu_2^{(1)}(f,g) + \Phi^{(1)}(f)g + f\Phi^{(1)}(g) 
        \label{eq:catequiv3.1}
\end{eqnarray}

That is, two first order deformations defined by $\mu_1^{(1)}$ and
$\mu_2^{(1)}$ are equivalent exactly when the $\mu_i^{(1)}$'s are
cohomologous. 

So we establish the folk theorem generalizing the classical result of
Gerstenhaber:

\begin{theorem} \label{1catclass}
The first order deformations of a $k$-linear category $\cal C$ are
classified up to equivalence by ${\frak H}^2({\cal C})$, the second
Hochschild
cohomology of the category.
\end{theorem}

Similarly, equations \ref{eq:functdef1}-\ref{eq:functdef3} become

\begin{eqnarray}
F^{(0)} & = & F \label{eq:functdef1.1}  \\
F^{(1)}(1_x) & = & 0  \,\, \makebox{\rm for all $x \in Ob({\cal C})$}
\label{eq:functdef2.1} \\
\lefteqn{F(\mu^{(1)}(f,g)) + F^{(1)}(fg) =} \nonumber \\
& & 
\nu^{(1)}(F(f), F(g)) + F^{(1)}(f)F(g) + F(f)F^{(1)}(g)
\label{eq:functdef3.1}
\end{eqnarray}
 
Again, as expected, in the case where the deformations of the source
and 
target categories are both trivial (i.e. $\mu^{(1)} \equiv 0$ and
$\nu^{(1)}
\equiv 0$), equation \ref{eq:functdef3.1} says that $F^{(1)}$ is a 
Hochschild 1-cocycle, while in the general case it cobounds
$F_*(\mu^{(1)}) - F^*(\nu^{(1)})$.  Or, put another way, 
$(\mu^{(1)}, \nu^{(1)}, F^{(1)})$ are a 1-cocycle in the cone on
$F^*(p_2) - F_*(p_1)$.

Weak equivalence between two first order deformations $(\hat{\cal
C}_i,
\hat{\cal D}_i,\hat{F}_i)$, $i = 1,2$, of $F:{\cal C}\rightarrow
{\cal D}$, 
then consists of equivalences of first order deformations of 
categories $\Gamma = Id_{\cal C} + \Gamma^{(1)}\epsilon$ from
$\hat{\cal C}_1$
to $\hat{\cal C}_2$, and  
$\Delta = Id_{\cal D} + \Delta^{(1)}\epsilon$ from $\hat{\cal D}_1$
to $\hat{\cal D}_2$, together with $\phi$ given by 
$\phi_x = 1_x + \phi^{(1)}_x$ satisfying 

\begin{eqnarray}
\Delta^{(1)}(F(f)) + F^{(1)}_1(f) + F(f)\phi^{(1)}_y & = &
     \phi^{(1)}_xF(f) + F^{(1)}_2(f) + F(\Gamma^{(1)}(f))
     \label{eq:functequiv2.1} 
\end{eqnarray}

\noindent for all $f:x\rightarrow y \in Arr({\cal C}).$

It is easy to see in the case of a strong equivalence, where
$\Delta^{(1)}$ and
$\Gamma^{(1)}$ are both zero, that this says $\phi$ cobounds the
difference
of $F^{(1)}_1$ and $F^{(1)}_2$.  Recalling a result from Gerstenhaber
and Schack \cite{GSdiag} makes it obvious what is happening in the
general case:  $(\Gamma^{(1)}, \Delta^{(1)}, \phi^{(1)})$ cobounds the
difference  
$(\mu^{(1)}_2, \nu^{(1)}_2, F^{(1)}_{2}) - (\mu^{(1)}_1, \nu^{(1)}_1,
F^{(1)}_{1})$
in the cone on $F^*(p_2) - F_*(p_1)$.  So we have

\begin{theorem} \label{1functclass}
The first order deformations of a $k$-linear functor 
$F:{\cal C}\rightarrow {\cal D}$ are classified up to weak equivalence
by the first cohomology ${\frak H}^1(F)$ of the deformation complex
of $F$.
\end{theorem}

Again this is the obvious generalization of the classical result for
algebras (or rather algebra homomorphisms).

For deformations $(\hat{\cal C}, \hat{\cal D}, \hat{F}, \hat{G},
\hat{\sigma})$
of a natural transformation $\sigma:F\Rightarrow G$ ,
for $F,G:{\cal C}\rightarrow {\cal D}$, equation \ref{eq:nattransdef}
becomes

\begin{eqnarray}
\lefteqn{\nu^{(1)}(F(f),\sigma_y) + F^{(1)}(f)\sigma_y +
F(f)\sigma^{(1)}_y  
=} \nonumber \\
 & &
     \nu^{(1)}(\sigma_x,G(f)) + \sigma^{(1)}_xG(f) +
\sigma_xG^{(1)}(f)
\label{eq:nattransdef.1}
\end{eqnarray} 

Here the cohomological interpretation of this equation is not
so clear.  However, once the somewhat baroque definition of the 
cochain map $\sigma \ddagger$ is recalled, it is easy to see that
this equation is the additional requirement on the fifth coordinate
to ensure that 
$(\mu^{(1)},\nu^{(1)},F^{(1)},G^{(1)},\sigma^{(1)})$ 
be a 0-cocycle
in the cone on $\sigma \ddagger$.

Similarly, equation \ref{eq:nattranswe} reduces to

\begin{eqnarray}
\phi^{(1)}_x\sigma_x + \sigma^{(1)}_{2 x} & 
     = & \sigma_x \psi^{(1)}_x + \sigma^{(1)}_{1 x}
\label{eq:nattranswe.1}
\end{eqnarray}

This is the condition on the $C^0(F,G)$ coordinate for 
$(\Gamma^{(1)},\Delta^{(1)},\psi^{(1)},\phi^{(1)},0)$ to cobound the
difference of the 
$(\mu^{(1)}_i,\nu^{(1)}_i,G^{(1)}_i,F^{(1)}_i,\sigma^{(1)}_i)$    
($i = 1,2$).  The other conditions requiring the other coordinates to
give equivalences of the deformations source and target functors and
categories give the remaining conditions, so we have

\begin{theorem} \label{1natclass}
The first order deformations of $\sigma:F\rightarrow G$ are
classified 
up to weak equivalence by the zeroth cohomology ${\frak H}^0(\sigma)$
of the
deformation complex of $\sigma$.
\end{theorem}

\section{Higher order deformations and obstructions without pasting}

Equations \ref{eq:catdef4}, \ref{eq:functdef3}, and
\ref{eq:nattransdef} are
the crucial defining conditions on deformations of categories,
functors and
natural transformations, respectively.  
The other equations are either normalization conditions 
to ensure correct behavior on identity arrows, conditions on the  
zeroth order term to ensure reduction to what is being deformed, or
conditions inherited from the deformation of a source or target.

In each case, we can solve the equation of index $n$ to separate
the terms involving coefficients of index $n$ from the other terms.

Equation \ref{eq:catdef4} gives

\begin{eqnarray} \label{compobstruction}
\lefteqn{\sum_{i+j=n \; 0\leq i,j < n} 
             \mu^{(i)}(a,\mu^{(j)}(b,c)) -
\mu^{(i)}(\mu^{(j)}(a,b),c) } 
           \nonumber \\
& =  & a\mu^{(n)}(b,c) - \mu^{(n)}(ab,c) + \mu^{(n)}(a,bc) -
\mu^{(n)}(a,b)c \\
& = & \delta \mu^{(n)}(a,b,c)
\end{eqnarray}

As expected, the left-hand side is the same formula as Gerstenhaber's 
obstructions for the deformation theory of associative algebras, and
the
condition that the next term satisfies is that it cobounds the
obstruction.
Letting $\langle \epsilon \rangle$ be the maximal ideal of
$k[\epsilon]/\langle\epsilon^{n}\rangle$ or $k[[\epsilon]]$ and
$\otimes$ denote $\otimes_{k}$ or its $\langle \epsilon \rangle$-adic
completion, as appropriate, these conditions for all $n$ can be neatly packaged into the single
condition that 

\[ \bar{\mu} 
= \sum \mu^{(i)}\epsilon^i \in C^2({\cal C})\otimes
\langle \epsilon \rangle \]

\noindent satisfies the  Maurer-Cartan equation

\[ \delta(\bar{\mu}) = \frac{1}{2}[\bar{\mu},\bar{\mu}], \]

\noindent where the (graded) Lie bracket is the obvious
generalization of
that given for algebras by Gerstenhaber,  and $\delta$ is the 
extension of the coboundary operator by linearity (and continuity),
or equivalently
the equation

\[ \delta(\bar{\mu}) = \bar{\mu}\{\bar{\mu}\}. \]

The presence of sources and targets is no impediment to the same
proof 
as in \cite{G1} that
each obstruction is always a 3-cocycle.

Similarly solving equation \ref{eq:functdef3} to separate index $n$
terms
gives

\begin{eqnarray} \label{functobstruction}
\lefteqn{\sum_{i+j=n, \, 0\leq i,j < n} F^{(i)}(\mu^{(j)}(f,g)) -
\sum_{k+l+m=n \, 0\leq k,l,m < n} \nu^{(k)}(F^{(l)}(f), F^{(m)}(g)) }
\nonumber \\
&  = & F(f)F^{(n)}(g) - F^{(n)}(fg) + F^{(n)}(f)F(g) \nonumber \\
& & \mbox{} - F(\mu^{(n)}(f,g) + \nu^{(n)}(F(f),F(g))
\end{eqnarray}

The right-hand side, is, of course, the $C^2(F)$ summand of
$\delta(\mu^{(n)},\nu^{(n)},F^{(n)})(f,g)$ in ${\frak C}^2(F)$. 

Letting $\bar{\mu}$ be defined as above, and $\bar{\nu}$ 
analogously, with $\bar{F} := \sum F^{(i)}\epsilon^i 
\in C^1(F)\otimes \langle \epsilon \rangle$, 
allows us collect all of the equations into the condition

\[
\delta(\bar{F},\bar{\mu},\bar{\nu}) 
= (\bar{F}\{\bar{\mu}\} - \bar{\nu}\{\bar{F}\}
-\bar{F}\bar{F} - \bar{\nu}\{\bar{F},\bar{F}\},
\bar{\mu}\{\bar{\mu}\}, \bar{\nu}\{\bar{\nu}\}), \]

\noindent where

\[ \delta(\bar{F},\bar{\mu},\bar{\nu}) 
= (\delta(\bar{F}) + F(\bar{\mu}) - \bar{\nu}(F,F), 
\delta(\bar{\mu}), \delta(\bar{\nu}) ) \]

\noindent since we are in the cone on $F^*(p_2) - F_*(p_1)$.  Notice
that since, $\bar{\mu}\{\bar{\mu}\} = \frac{1}{2}[\bar{\mu},
\bar{\mu}]$,
and similarly for $\bar{\nu}$, the last two coordinates are simple
restating the Maurer-Cartan equation for the deformation of the source
and target categories.   

Finally, solving equation \ref{eq:nattransdef} to separate index $n$
terms gives

\begin{eqnarray} \label{natobstruction}
\lefteqn{\sum_{i+j+k=n \, 0 \leq i,j,j < n} 
\nu^{(i)}(F^{(j)}(f),\sigma^{(k)}_y)
- \nu^{(i)}(\sigma^{(j)}_x,G^{(k)}(f))} \nonumber \\
& = & \nu^{(n)}(\sigma, G(f)) - \nu^{(n)}(F(f), \sigma) \nonumber \\
& & + \sigma^{(n)}G(f) - F(f)\sigma^{(n)} 
+ \sigma G^{(n)}(f) - F^{(n)}(f)\sigma
\end{eqnarray}

Here it is easy to verify that the right-hand side is the
$C^1(F,G)$ summand of 

\[ \delta(\mu^{(n)}, \nu^{(n)}, F^{(n)}, G^{(n)}, \sigma^{(n)}) \]
 
\noindent in
${\frak C}^1(\sigma)$. 

Again all instances of this equation can be collected into
a single equation relating the differential and the cup-like and
brace-like 
compositions:  letting $\bar{\mu}, \bar{\nu}$, and $\bar{F}$ be as
above, and $\bar{G}$ be defined similarly, and letting
$\bar{\sigma} = \sigma^{(i)}\epsilon^i$,
this equation becomes

\[ p_5(\delta(\bar{\sigma},\bar{F},\bar{G},\bar{\mu},\bar{\nu}))
= \delta(\bar{\sigma}) - \bar{\nu}\{\bar{F},\sigma\} -
\bar{\nu}\{\sigma,\bar{G}\}
-\bar{\nu}\{\bar{\sigma}\} + \bar{F}\bar{\sigma} - \bar{\sigma}\bar{G}
-\bar{\nu}\{\bar{F},\bar{\sigma}\} - \bar{\nu}\{\bar{\sigma},\bar{G}\}
\] 
 
\noindent The conditions on the other coordinates are those given
previously. 

Equations \ref{compobstruction}, \ref{functobstruction} and
\ref{natobstruction}, then identify the obstructions to extending an
$n-1^{st}$ order deformation to an $n^{th}$ order
deformation of a category, functor or natural transformation,
respectively.

As expected, we have

\begin{theorem} \label{obstructionsareclosed} The obstruction

\[ \omega_{\cal A}^n := \sum_{i+j=n \; 0\leq i,j < n} 
             \mu^{(i)}(a,\mu^{(j)}(b,c)) -
\mu^{(i)}(\mu^{(j)}(a,b),c)  \]

\noindent (resp.

\[ \omega_{F}^{n} : = (\omega_{A}^{n}, \omega_{B}^{n}, 
\upsilon_{F}^{n})\]

\noindent where the first and second coordinates are obstructions to 
deforming the categories $\cal A$ and $\cal B$ respectively, and

\[ \upsilon_{F}^{n} := 
\sum_{i+j=n, \, 0\leq i,j < n} F^{(i)}(\mu^{(j)}(f,g)) -
\sum_{k+l+m=n \, 0\leq k,l,m < n} \nu^{(k)}(F^{(l)}(f), F^{(m)}(g)), 
\]

	
\[ \omega_{\sigma}^{n} : = (\omega_{A}^{n}, \omega_{B}^{n}, 
\upsilon_{F}^{n}, \upsilon_{G}^{n}, \upsilon_{\sigma}^{n})\]

\noindent where the first four coordinates are as described above and

\[ \upsilon_{\sigma}^{n} := \sum_{i+j+k=n \, 0 \leq i,j,j < n} 
\nu^{(i)}(F^{(j)}(f),\sigma^{(k)}_y)
- \nu^{(i)}(\sigma^{(j)}_x,G^{(k)}(f)) \; \mbox{\rm{\em )}}\] 


\noindent is closed in ${\frak C}^\bullet({\cal A})$ (resp.
${\frak C}^\bullet(F)$, ${\frak C}^\bullet(\sigma)$).
\end{theorem} 

\noindent{\sc Proof.}  The statements concerning the obstructions
in the case of categories (resp. functors) follows by the same proof
given by Gerstenhaber \cite{G1,G2} in the case of algebras 
(resp. Gerstenhaber and Schack \cite{GSdiag} in the case of algebra
homomorphisms), as does the vanishing of all but the last coordinate
in $\delta(\omega_{\sigma}^n)$.

It thus remains only to show that the last coordinate of 
$\delta(\omega_{\cal A}^n,\omega_{\cal
B}^n,\omega_G^n,\omega_F^n,\omega_{\sigma}^n)$ vanishes.

The coboundary is

\begin{eqnarray*}
\lefteqn{ \delta_{F,G}(\upsilon_{\sigma}^n) - \sum_{i+j=n} \nu^{(i)}
\{\nu^{(j)}\} \{\sigma\} - } \\
& &
\sigma_x \left[ \sum_{i+j=n} G^{(i)}(\mu^{(j)}(f,g)) - \sum_{i+j+k=n}
\nu^{(i)}(G^{(j)}(f),G^{(k)}(g)) \right]
+ \\
& & \left[ \sum_{i+j=n} F^{(i)}(\mu^{(j)}(f,g)) - \sum_{i+j+k=n}
\nu^{(i)}(F^{(j)}(f),F^{(k)}(g)) \right] \sigma_z  
\end{eqnarray*}

The key to the straightforward but tedious calculation which shows
this vanishes is to immediately rewrite $\sum_{i+j=n}
G^{(i)}(\mu^{(j)}(f,g))$ using equation \ref{eq:nattransdef} (or
equivalently \ref{natobstruction}), cancel the terms involving
$F^{(i)}$'s with those already in the original expression, then
rewrite
the terms still involving $\mu^{(j)}$'s using equation
\ref{eq:functdef3}
(or \ref{functobstruction}).   

Using equations \ref{compobstruction}, \ref{functobstruction}, and
\ref{natobstruction} the number of terms can be steadily reduced
(though at four points in the calculation as carried out by the
author, \ref{natobstruction} must be used to replace three sums by
three others).  At two points the naturality of
$\sigma$ must be used. $\blacksquare$

The reader intent on recovering the complete calculation for him- or
herself is advised to first carry out the case of $n=2$, where
equations \ref{compobstruction}, \ref{functobstruction}, and
\ref{natobstruction} are simply cocycle conditions.  $\blacksquare$

\section{Deformations induced by single compositions}

Since deformations of categories, functors and natural transformations
are themselves categories, functors and natural transformations, it is
clear that the usual operations in the 2-category $k[[\epsilon]]-cat$
or
$k[\epsilon]/\langle \epsilon^2 \rangle-cat$ induce operations on the 
deformations, and in particular, that deformations of the parts
of a composable pasting diagram induce a deformation of its pasting
composition.

Relatively trivial calculations establish formulas for the induced
deformation in the case of individual compositions:

\begin{proposition} \label{1compositionofdeformations}
If ${\cal A}\stackrel{F}{\rightarrow} {\cal B}
\stackrel{G}{\rightarrow} 
{\cal C}$ is a composition of $k$-linear functors with composite
$\Phi:{\cal A}\rightarrow {\cal C}$, and $\tilde{\cal A},
\tilde{\cal B}, \tilde{\cal C}, \tilde{F},$ and $\tilde{G}$ are
deformations
of its parts, with $\tilde{F} = \sum F^{(i)}\epsilon^i$ and 
$\tilde{G} = \sum G^{(i)}\epsilon^i$, then 

\[ \Phi^{(i)} := \sum_{j+k = i} G^{(j)}(F^{(k)}) \]

\noindent are the terms of a deformation of $\Phi$, called {\em the
induced
deformation of $\Phi$}.
\end{proposition}

\begin{proposition}
If $F:{\cal A}\rightarrow {\cal B}$ is a $k$-linear functor, and
$\sigma:G \Rightarrow H$ is a natural transformation between
$k$-linear functors $G,H:{\cal B}\rightarrow {\cal C}$, and
$\tilde{\cal A}, \tilde{\cal B}, \tilde{\cal C}, \tilde{F}, 
\tilde{G}, \tilde{H},$ and $\tilde{\sigma}$ are deformations of
its parts, with $\tilde{\sigma}_X = \sum \sigma^{(i)}_X \epsilon^i$,
then $\tilde{\tau} = \sum \sigma_{F(X)}^{(i)}\epsilon^{i}$ 
is a deformation of $\tau = \sigma_F$ as
a natural transformation from the induced deformation of $G(F)$ to 
the induced deformation of $H(F)$. 
\end{proposition}

\begin{proposition}
If $F,G:{\cal A}\rightarrow {\cal B}$ are $k$-linear functors, and
$\sigma:F \Rightarrow G$ is a natural transformation between them, and
$H:{\cal B}\rightarrow {\cal C}$ is a $k$-linear functor and
$\tilde{\cal A}, \tilde{\cal B}, \tilde{\cal C}, \tilde{F}, 
\tilde{G}, \tilde{H},$ and $\tilde{\sigma}$ are deformations of
its parts, with $\tilde{\sigma}_X = \sum \sigma^{(i)}_X \epsilon^i$,
and $\tilde{F} = \sum F^{(i)} \epsilon^i$, then $\tilde{\tau} = \sum
\tau^{(i)} \epsilon^i$
given by

\[ \tau^{(i)}_X := \sum_{j+k=i} H^{(j)}(\sigma^{(k)}_X) \]

\noindent is a deformation of $H(\sigma)$ as a natural transformation
from the
induced deformation of $H(F)$ to the induced deformation of $H(G)$.
\end{proposition}

Finally, we have

\begin{proposition} \label{2compositionofdeformations}
If $F,G,H:{\cal A}\rightarrow {\cal B}$ are $k$-linear functors, and
$\sigma:F \Rightarrow G$ and $\tau:G \Rightarrow H$ are natural 
transformations with composite $\phi:F\Rightarrow G$, 
and $\tilde{\cal A}, \tilde{\cal B}, \tilde{F}, 
\tilde{G}, \tilde{H},\tilde{\sigma},$ and $\tilde{\tau}$ are 
deformations of the parts, then 

\[ \phi^{(i)} := \sum_{k+j+l = i} \nu^{(k)}(\sigma^{(j)},\tau^{(l)}) 
\epsilon^i \]

\noindent defines a deformation of the composite $\phi$.
\end{proposition}

Observe that in each proposition, the terms of the induced
deformation can be expressed
in terms of the brace-like 1-composition, being $G^{(j)}\{F^{(k)}\}$,
$\sigma^{(i)}\{F\}$, $F^{(j)}\{\sigma^{(k)}\}$ and
$\nu^{(k)}\{\sigma^{(j)}, \tau^{(l)}\}$, respectively. 

It follows from Power's Pasting Theorem that deformations of all
parts of a composable pasting diagram induce a deformation of the
composite.

\section{The cohomology of $k$-linear pasting diagrams}

What remains now is to fit the parts introduced thus far together.
To do this we must return to the description of a pasting diagram.

First, we should note that not just the image of the diagram, the
categories, functors and natural transformations involved, but the
`shape' of the diagram will matter a great
deal:  even if the labeling of the pasting scheme includes
coincidences, the different copies of the same category, functor, or 
natural transformation may be deformed independently.

Thus, for instance, in the diagram

\begin{center}
\setlength{\unitlength}{4144sp}%
\begingroup\makeatletter\ifx\SetFigFont\undefined%
\gdef\SetFigFont#1#2#3#4#5{%
  \reset@font\fontsize{#1}{#2pt}%
  \fontfamily{#3}\fontseries{#4}\fontshape{#5}%
  \selectfont}%
\fi\endgroup%
\begin{picture}(2055,930)(1831,-1231)
\thinlines
{\color[rgb]{0,0,0}\put(2116,-646){\vector( 1, 0){1620}}
}%
{\color[rgb]{0,0,0}\put(2116,-961){\vector( 1, 0){1620}}
}%
\put(2791,-466){\makebox(0,0)[lb]{\smash{{\SetFigFont{12}{14.4}{\rmdefault}{\mddefault}{\updefault}{\color[rgb]{0,0,0}$F$}%
}}}}
\put(1846,-871){\makebox(0,0)[lb]{\smash{{\SetFigFont{12}{14.4}{\rmdefault}{\mddefault}{\updefault}{\color[rgb]{0,0,0}$\cal
A$}%
}}}}
\put(3871,-856){\makebox(0,0)[lb]{\smash{{\SetFigFont{12}{14.4}{\rmdefault}{\mddefault}{\updefault}{\color[rgb]{0,0,0}$\cal
A$}%
}}}}
\put(2776,-1216){\makebox(0,0)[lb]{\smash{{\SetFigFont{12}{14.4}{\rmdefault}{\mddefault}{\updefault}{\color[rgb]{0,0,0}$F$}%
}}}}
\end{picture}%
\end{center}

\noindent a deformation will involve two, not necessarily equal,
deformations of the category $\cal A$, and two, not necessarily
equal, deformations of the functor $F$.  

If we were interested, instead, in a deformation of the endofunctor
$F$, as an endofunctor, that is a single deformation of $\cal A$, and
a deformation of $F$ in which both source and target are deformed by
this deformation, we would, instead be studying the deformations of
the diagram

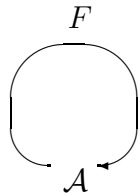
\begin{center}
\setlength{\unitlength}{4144sp}%
\begingroup\makeatletter\ifx\SetFigFont\undefined%
\gdef\SetFigFont#1#2#3#4#5{%
  \reset@font\fontsize{#1}{#2pt}%
  \fontfamily{#3}\fontseries{#4}\fontshape{#5}%
  \selectfont}%
\fi\endgroup%
\begin{picture}(772,1155)(2217,-1096)
\thinlines
{\color[rgb]{0,0,0}\put(2746,-931){\vector(-1, 0){0}}
\put(2746,-582){\oval(470,698)[br]}
\put(2604,-582){\oval(754,754)[tr]}
\put(2604,-582){\oval(756,754)[tl]}
\put(2461,-582){\oval(470,698)[bl]}
}%
\put(2536,-1081){\makebox(0,0)[lb]{\smash{{\SetFigFont{12}{14.4}{\rmdefault}{\mddefault}{\updefault}{\color[rgb]{0,0,0}$\cal
A$}%
}}}}
\put(2566,-106){\makebox(0,0)[lb]{\smash{{\SetFigFont{12}{14.4}{\rmdefault}{\mddefault}{\updefault}{\color[rgb]{0,0,0}$F$}%
}}}}
\end{picture}%
\end{center}

In light of this discussion, it should be clear that the relevant
cochain complex for the deformation of an entire pasting diagram will
arise by iterated mapping cone constructions from cochain complexes
associated to the
labels of the various cells of the diagram.

In particular, the groups of cochains for a diagram 

\[D:G\rightarrow 3-computad(k-Cat), \]

\noindent are given by

\begin{eqnarray*}
 {\frak C}^\bullet(D) & = & \bigoplus_{v \in G_0} C^{\bullet+3}(D(v))
\oplus
\bigoplus_{e \in G_1} C^{\bullet+2}(D(e)) \oplus \\
& & \bigoplus_{f \in G_2}
C^{\bullet+1}(dom(D(f)), cod(D(f)) \oplus\\
& & \bigoplus_{s \in G_3} C^{\bullet}(dom(cod(D(s))),
cod(cod(D(s)))) \end{eqnarray*}

\noindent Here, of course each $D(v)$ is a linear category, each
$D(e)$ is a functor, as are $cod(D(f))$, $dom(D(f))$, $cod(cod(D(s)))
=
cod(dom(D(s)))$, and $dom(cod(D(s))) = dom(dom(D(s)))$, these last
four being the composition of the $D(e)$'s along the codomain of 
a 2-arrow (resp. the domain of a 2-arrow, the common codomain of the
domain and codomain of a 3-cell, and the common domain of the domain
and codomain of a 3-cell).  (Recall 3-cells have composable pasting
diagrams as domain and
codomain, and merely assert the equality of the composites.)  

Notice
here we are abusing notation somewhat by not explicitly applying 
$comp(-)$ to the domains and codomains of 2- and twice iterated 
domains and codomains of 3-cells, which 
properly are composable 1-pasting schemes.  This
should cause no confusion, as the Hochschild complexes are only 
defined for pairs of parallel 1-arrows, not for composable 1-pasting
schemes.  We retain this convention throughout what follows.

The tricky thing is to succinctly describe the coboundary maps.  The
dimension shifts hint at the construction:  the coboundary maps will
arise from an iterated mapping cone construction.

First we need 

\begin{proposition} \label{gencochain1comp}
If $D$ is a composable 1-pasting diagram, and $F$ is a
1-arrow therein, then the pre- and post-composition cochain maps of 
Proposition \ref{cochain1comp} induce a unique cochain map

\[ \wp^F_D: C^\bullet(F)\rightarrow C^\bullet(comp(D)), \]

\noindent where $comp(D)$ denotes the composition of the arrows in
$D$.

Similarly if $D$ is a 1-pasting diagram which is the union of two
composable 1-pasting diagrams $D_F$ and $D_G$, which are identical
except for the label on one element of the underlying $G_1$, which
is $F$ and $G$ respectively, then the pre- and post-composition 
cochain maps of 
Proposition \ref{cochain1comp} induce a unique cochain map

\[ \wp^{F,G}_D: C^\bullet(F,G)\rightarrow
C^\bullet(comp(D_F),comp(D_G)). \]
\end{proposition}

\noindent{\sc Proof.}
The map is constructed simply by iterated application of the cochain 
maps of Proposition \ref{cochain1comp}.  It is unique by associativity
of 1-composition. $\blacksquare$

Similarly, with the quite non-trivial proof of uniqueness in
Power's Pasting Theorem \cite{Pow} replacing
the rather trivial proof that associativity implies uniqueness of
iterated compositions, we have

\begin{proposition} \label{gencochain2comp}
If $D$ is a composable 2-pasting diagram, and $\sigma$ is a 
2-arrow therein, the pre- and post-composition cochain maps of
Propositions \ref{cochain1comp} and \ref{cochain2comp} induce
a unique cochain map

\[ \wp^\sigma_D : C^\bullet(comp(dom(\sigma)),comp(cod(\sigma)))
\rightarrow C^\bullet(comp(dom(D)),comp(cod(D))) \]
\end{proposition}

Armed with these results, we can now proceed to construct the
coboundaries:

First, observe that there is a cochain map

\begin{eqnarray*}
\lefteqn{ \kappa_1 : =  \sum_{e \in G_1} \iota_e( D(e)^*( p_{cod(e)}
) ) -
\sum_{e \in G_1} \iota_e( D(e)_*( p_{dom(e)} ) ):} \\ 
& &
\bigoplus_{v \in G_0} C^\bullet( D(v) )\rightarrow
\bigoplus_{e \in G_1} C^\bullet( D(e) ). \end{eqnarray*}

The components of this map, are, of course, the cochain maps arising
in
the construction of ${\frak C}^\bullet ( D(e) )$ as in Section
\ref{nopasting}, and the fact that it is a cochain map follows from
this.

The cone on this map, which we denote ${\frak C}_1^\bullet(D)$ has
cochain groups

\[
\bigoplus_{v \in G_0} C^{\bullet+1}( D(v) ) \oplus
\bigoplus_{e \in G_1} C^{\bullet}( D(e) ) 
\]

This construction is simply a replication for each 1-cell of the
pasting
diagram of that given in Section \ref{nopasting}.   1-cocycles in 
${\frak C}_1^\bullet(D)$ classify simultaneous deformation of all
objects and functors in the diagram (recalling our earlier warning 
that occurences of a functor or category at a different place in 
the diagram may be deformed differently.

For 2-cells, we
cannot simply replicate the earlier construction, because in general
the source and target are composable 1-pasting diagrams, rather than
1-cells.  We therefore proceed in two steps, the second of which
corresponds
to replicating the construction in Section \ref{nopasting}, while the
first
involves the maps $\wp^F_D$.

We then have

\begin{proposition}  For any pasting diagram $G$,

\[
    \wp : = 
    \sum_{v \in G_{0}} i_{v}(p_{v}) +
    \sum_{e \in \Delta \in diag_1(G)} 
i_{\Delta}(\wp^{D(e)}_{D(\Delta)}(p_{e})):  {\frak C}_1^\bullet(D) 
\rightarrow  \]
\[ \bigoplus_{v \in G_{0}} C^{\bullet+1}(D(v)) \oplus
\bigoplus_{\Delta \in diag_{1}(G)} C^\bullet (comp(D(\Delta))),
\]

\noindent  is a cochain map, 
where $diag_1(G)$ denotes the set of composable 1-pasting
diagrams in $G$.  It assigns to a 1-cocycle in ${\frak C}_1^\bullet(D)$ 
which names a deformation of all categories and functors in the 
diagram, a 1-cocycle in the target complex whose $\Delta$ summand
gives the deformations of the domain and codomain categories of
$D(\Delta)$ and the deformation of the $comp(D(\Delta))$ induce by
the deformations of the composed functors.
\end{proposition}

\noindent{\sc Proof.} That the first two coordinates of each summand
behave correctly is immediate.  That the third coordinate in each summand 
behaves correctly follows from the fact that the maps of Proposition
\ref{gencochain1comp} are cochain maps, and easily verified the cancellation of
terms involving images of cochains associated to an intermediate 
category in the composable pasting diagram $\Delta$.
$\blacksquare$\smallskip

And by replicating the construction of Section \ref{nopasting}, there
is a cochain map 

\begin{eqnarray*} \lefteqn{ \kappa_2 := \sum_{\sigma \in G_2}
i_\sigma( D(\sigma)\ddagger(\pi_\sigma))
:} \\ & &
\bigoplus_{v \in G_{0}} C^{\bullet+1}(D(v)) \oplus
\bigoplus_{\Delta \in diag_{1}(G)} C^\bullet (comp(D(\Delta)))
\rightarrow \\ & &
\bigoplus_{\sigma \in G_2} C^\bullet(dom(D(\sigma))),
cod(D(\sigma)))), \end{eqnarray*}

\noindent where $\pi_\sigma$ is the projection from the source onto
its 
summands constituting 

\[ C^{\bullet}(dom(cod(D(\sigma)))
\stackrel{dom(D(\sigma))}{cod(D(\sigma))} cod(cod(D(\sigma)))). \]


The cone on the composite $\kappa_2(\wp)$, which we denote
${\frak C}_2^\bullet(D)$ then has cochain groups

\[\bigoplus_{v \in G_0} C^{\bullet+2}(D(v)) \oplus
\bigoplus_{e \in G_1} C^{\bullet+1}(D(e)) \oplus \\
 \bigoplus_{f \in G_2}
C^{\bullet}(dom(D(f)), cod(D(f)).\]

We now need one last cone construction to allow the presence of
3-cells
in a pasting diagram to enforce equality between its source and
target:

Define a cochain map

\begin{eqnarray*} \lefteqn{\kappa_3 := \sum_{c \in G_3, \; f\in cod(c)}
\iota_{c}(\wp_{cod(c)}^{D(f)}(p_{f})) - \sum_{c \in G_3, \; f\in dom(c)}
\iota_{c}(\wp_{dom(c)}^{D(f)}(p_{f})):} \\ & & 
{\frak C}_2^\bullet(D) \rightarrow
\bigoplus_{c \in G_3} C^\bullet(dom(cod(D(s))), cod(cod(D(s)))). 
\end{eqnarray*}

The cone on this map is then the desired cochain complex, which
we denote ${\frak C}^\bullet(D)$, and call {\em the deformation 
complex of the pasting scheme $D$}.  

Observe that it is this last 
step that obliged us to define ${\frak C}^{\bullet}(1_{\sigma})$ to 
be a cone on the difference of two projections, rather than simply 
${\frak C}^{\bullet}(\sigma)$:  in the context of a pasting scheme, 
the source and target of a 3-cell are, in general, pasting 
compositions of natural transformations, which must each be deformed.
The only convenient artifice within the cohomological framework for
enforcing equality of the induced deformations on the two sides of a
3-cell is to combine the maps $\wp^{\sigma}_{D}$ of Proposition
\ref{gencochain2comp} with the last cone construction of Section
\ref{Hochschild}.

We then have

\begin{theorem} \label{classificationfordiagrams}

First order deformations of the pasting diagram $D$ are classified up
to equivalence by ${\frak H}^{-1}(D)$, the negative-first
cohomology of the deformation complex ${\frak C}^\bullet(D)$.

\end{theorem}

\noindent and

\begin{theorem} \label{obstructionsfordiagrams}

Given an $n-1^{st}$-order deformation of a pasting diagram $D$, there
is a cocycle in ${\frak C}^0(D)$, each direct summand of which
is given by the formula for the obstruction to deforming the
label on the cell of the computad indexing the direct summand given
in Theorem \ref{obstructionsareclosed}.  This cocycle is
the obstruction to extending the deformation to an $n^{th}$ order
deformation, and if it vanishes in cohomology, any 0-cochain
cobounding it gives the degree $n$ term of an $n^{th}$ order
deformation extending the given deformation.

\end{theorem}

Notice that the somewhat strange cohomological dimensions are
correct:  the cohomological dimension in the cone corresponds to the
cohomological dimension in the groups associated to the 3-cells of
the 3-computad, so the $-1$-cocycles in the cone have coordinates
which are a $-1$-cochain for each 3-cell (necessarily 0, indicating
equality between its source and target), and, as expected, a
$0$-cochain for each 2-cell, a
$1$-cochain on each 1-cell, and a $2$-cochain on each 0-cell,
collectively statisfying the cocycle condition for the iterated
mapping cone.

\noindent {\sc Proofs:} Both results would follow immediately from
Theorems \ref{1catclass}, \ref{1functclass}, \ref{1natclass} and
\ref{obstructionsareclosed}  and Propositions
\ref{1compositionofdeformations} through
\ref{2compositionofdeformations}, (albeit with the cohomological
dimension shifted up 1) were it not for the presence of 3-cells
indicating commutative parts of the pasting diagram.  However, the
presence of 3-cells, their own labels being undeformable, 
enforces the equality of the induced deformations on their source and
target, since a 0-cocycle (in this case the difference between either 
the induced first-order deformations of the source and target, or the 
difference of the cocycles induced by the obstruction on each face of 
the source and targed is trivial in cohomology if and only if it is 
zero, there being only the zero -1-cochain to cobound it.
$\blacksquare$

\section{Prospects}

As regards the first motivation for this paper: in work in progress,
a doctoral student under
the author is attempting the construction of a cohomology theory
governing the simultaneous deformation of the composition, arrow-part
of the monoidal product, and structure maps of a monoidal category.
It appears that the deformations are governed by the total complex of
a 'multicomplex'--a bigraded object which 'looks like a spectral
sequence with all the pages smashed together'.  The theory of the
present paper gives the (0,1)-differentials in one direction, while
the (1,0)-differentials are given by the differentials of
\cite{Y.book, Y.abkan}.

As regards the second motivation: it is easy to see that $k$-linear
stacks are a special case of pasting diagrams, indeed, $k$-linear
pre-stacks are more or less the same thing as $k$-linear pasting
diagrams.  It is clear that deformation of a pre-stack (as a pasting
diagram) cannot create new descent data, since any descent data in an
order $n$ deformation must be descent data at all lower orders, as
quotienting by powers of $\epsilon$ will preserve commutativity.  On
the other hand, in general deformation can destroy commutativity, and
thus descent data. It appears natural conjecture that
any effective descent data which is not destroyed by a given
deformation remains effective in the deformation.

The author plans to prove this conjecture in subsequent work.

\newpage

\begin{flushleft}
\bibliography{pastingdef}
\end{flushleft}

\end{document}